\documentclass[11pt]{article}
\usepackage{amsmath, amsfonts, amssymb, amsthm, a4wide}
\usepackage[arrow]{xy}
\numberwithin{equation}{section}
\title{Gauss Decomposition of the Yangian $Y(\mathfrak{gl}_{m|n})$}
\author{Lucy Gow\footnote{lucyg@maths.usyd.edu.au} \\ {\normalsize{University of Sydney}}}
\begin{document}
\maketitle
\newcommand{\pa}[1]{\overline{#1}}
\newcommand{\gl}{\mathfrak{gl}}
\newtheorem{lemma}{Lemma}[section]
\newtheorem{corollary}[lemma]{Corollary}
\newtheorem{proposition}[lemma]{Proposition}
\newtheorem{defn}{Definition}[section]
\newtheorem{theorem}{Theorem}
\newtheorem{conjecture}[theorem]{Conjecture}
\newtheorem*{dualform}{Dual Form of the Quantum Berezinian}
\newtheorem{remark}{Remark}[section]
\maketitle
We describe a Gauss decomposition for the Yangian $Y(\gl_{m|n})$ of the general linear Lie superalgebra.
This gives a connection between this Yangian and the Yangian of
the classical Lie superalgebra $Y(A({m-1, n-1}))$ (with $m\ne n$) defined and studied in papers by
Stukopin, and suggests natural definitions for the Yangians $Y(\mathfrak{sl}_{n|n})$ and $Y(A(n,n))$. 
We also show that the coefficients of the quantum Berezinian generate
the centre of the Yangian $Y(\gl_{m|n})$.  This was conjectured by Nazarov in 1991.

\section{Introduction}
The Yangian $Y (\mathfrak{gl}_{m|n})$
is the $\mathbb{Z}_2$-graded
associative algebra over $\mathbb{C}$ with generators 
\[
\{t_{ij}^{(r)}\,| \; 1\le i,j \le m+n; r\ge 1\}
\]
and defining relations
\begin{equation}
\label{dr4}
[t_{ij}^{(r)}, t_{kl}^{(s)}] = (-1)^{\pa{i}\,\pa{j} + \pa{i}\,\pa{k} + \pa{j}\,\pa{k}}
\sum_{p=0}^{\mathrm{min}(r,s) -1}(t_{kj}^{(p)} t_{il}^{(r+s-1-p)} - t_{kj}^{(r+s-1-p)}t_{il}^{(p)}).
\end{equation} where $\pa{i}$ is the parity of the index $i$. 
We take $\pa{i} = 0$ for $i\le m$; and $\pa{i}=1$ for $i \ge m+1$.
(We write square brackets for the super-commutator).
We define the formal power series
\begin{equation*}
t_{ij} (u) = \delta_{ij} + t_{ij}^{(1)} u^{-1} + t_{ij}^{(2)}u^{-2} + \ldots
\end{equation*}
and a matrix
\begin{equation}\label{T}
T(u) = \sum_{i,j =1}^{m+n} t_{ij} (u) \otimes E_{ij} \; (-1)^{\pa{j} (\pa{i} +1)}
\end{equation}
where $E_{ij}$ is the standard elementary matrix. (Here we identify an operator 
$ \sum A_{ij}\otimes E_{ij} \; (-1)^{\pa{j} (\pa{i} +1)}$ in 
$Y(\gl_{m|n})[[u^{-1}]]\otimes \mathrm{End} \, \mathbb{C}^{m|n}$ with the matrix $\left( A_{ij} \right)_{i,j
= 1}^{m+n}$. 
The extra sign ensures that the product of two matrices can still be calculated in the usual way).
Then, as for the Yangian $Y(\mathfrak{gl}_{n})$ (see for example \cite{BK,MNO}),
 the defining relations may be expressed by the matrix product
\begin{equation*}
R(u-v) T_1 (u) T_2 (v) = T_2 (v) T_1 (u) R(u-v)
\end{equation*}
where
\begin{equation*}
R(u-v) = 1 - \frac{1}{(u-v)} P_{12}
\end{equation*}and $P_{12}$ is the permutation matrix:
$P_{12} = \sum_{i,j =1}^{m+n} E_{ij} \otimes E_{ji} (-1)^{\pa{j}}.$
We also have the following equivalent form of the defining relations:
\begin{equation}\label{dr2}
[t_{ij}(u),\;t_{kl}(v)] \;=\; \frac{(-1)^{\pa{i}\,\pa{j}+\pa{i}\,\pa{k}+\pa{j}\,\pa{k}}}{(u-v)}(t_{kj}(u)
t_{il} (v) - t_{kj} (v) t_{il} (u) ).
\end{equation}
The Yangian $Y(\gl_{m|n})$ is a Hopf algebra with comultiplication
\begin{equation}
\Delta : t_{ij}(u)  \mapsto \sum_{k=1}^{m+n} t_{ik}(u)  \otimes t_{kj}(u),
\end{equation}
antipode $S: T(u) \mapsto T(u)^{-1}$ and counit $\epsilon : T(u) \mapsto 1$. Throughout this article we
observe the following notation for entries of the inverse of the matrix $T(u)$: 
\begin{equation*}
T(u)^{-1} =: \left(t'_{ij}(u) \right)_{i,j=1}^{n}.
\end{equation*}
A straightforward calculation yields the following relation in $Y(\mathfrak{gl}_{m|n})$:
{
\begin{eqnarray}\label{useful}
[t_{ij}(u), t'_{kl}(v)]
= \frac{(-1)^{\pa{i}\,\pa{j} +\pa{i}\,\pa{k} + \pa{j}\,\pa{k}}}{(u-v)}\cdot (\; \delta_{kj} \sum_{s=1}^{m+n}
t_{is}(u) t'_{sl} (v) -\;\delta_{il} \sum_{s=1}^{m+n} t'_{ks}(v)t_{sj}(u))  .
\end{eqnarray}
}
We may define two different filtrations on the Yangian $Y(\gl_{m|n})$.  These are defined by setting the degree of a
generator as follows:
\begin{equation}\label{filtrations}
\mathrm{deg}_1 (t_{ij}^{(r)}) = r; \quad \quad \mathrm{deg}_2 (t_{ij}^{(r)}) = r-1.
\end{equation} 
Let gr$_1 Y(\gl_{m|n})$ and gr$_2 Y(\gl_{m|n})$, respectively, denote the corresponding graded algebras.

There is an injective homomorphism $\iota: U(\gl_{m|n}) \to Y(\gl_{m|n})$ given by 
\[ \iota: E_{ij} \mapsto t_{ij}^{(1)} (-1)^{\pa{i}}.\]
The injectivity of $\iota$ follows from the fact that its composition with a surjective homomorphism
$\;\pi:Y(\gl_{m|n}) \to U(\gl_{m|n})$ is the identity map on $U(\gl_{m|n})$.  The map $\pi$ is given as
follows:
\begin{equation}
\label{projection}
\pi: t_{ij}(u) \mapsto \delta_{ij} + E_{ij} (-1)^{\pa{i}} u^{-1}
\end{equation}
Thus we regard the universal enveloping algebra $U(\gl_{m|n})$ as a subalgebra of $Y(\gl_{m|n})$.

The Yangian $Y(\gl_{m|n})$ was introduced in \cite{Nazarov}. It has applications in mathematical physics
because it describes
symmetry in integrable models of Calogero-Sutherland systems~\cite{AhnKoo, JWW},
 superstrings in $AdS_{5} \times S^{5}$~\cite{HatsudaYoshida}, and in the hierarchy of a form of the non-linear
super-Schr\"{o}dinger equation with $m$ bosons and $n$ fermions \cite{CaudrelierRagoucy}.  The centre of the
Yangian $Y(\gl_{m|n})$ is conveniently described using a formal power series called the quantum Berezinian
(see Section \ref{centre}).

Vladimir Stukopin \cite{StukopinA, Stukopin} has introduced Yangians for classical simple Lie
superalgebras.  In this article we provide a new presentation for the Yangian
$Y(\gl_{m|n})$ that allows us to relate it to the Yangian $Y(A({m-1, n-1}))$ (for $m\ne n$) studied by
Stukopin.
This leads us to introduce a natural definition of the Yangian $Y(\mathfrak{sl}_{n|n})$
as a subalgebra of the Yangian $Y(\gl_{n|n})$, as well as a definition of $Y(A(n-1, n-1))$ 
(see Section \ref{slmnsection}).  
The Yangian that features in $D=4$ superconformal Yang-Mills theory~\cite{DolanNappiWitten} is that
associated with the supergroup $PSU(4,4)$, which has a Lie superalgebra of type $A(3,3)$, so the results
presented here may be relevant.  

Our paper follows similar treatments of the Yangian $Y(\gl_{N})$ given in papers by
Brundan and Kleshchev, and Cramp\'{e}, and
of the super-Yangian $Y(\gl_{1|1})$ in the work of Jin-fang Cai, Guo-xing Ju,  Ke Wu and Shi-kun
Wang (see \cite{BK,CJWW,Crampe}). 
\section{The Poincar\'{e}-Birkhoff-Witt Theorem for Super Yangians}  
In this section we prove the Poincar\'{e}-Birkhoff-Witt theorem for the Yangian $Y(\gl_{m|n})$.
The proof is based very closely on that of the corresponding theorem for $Y(\gl_N)$ given in \cite{BK}.

For each positive integer $l\ge 1$, we define a homomorphism 
\begin{displaymath}
\kappa_l := (\pi \otimes \cdots \otimes \pi) \circ \Delta^{(l)} : Y(\gl_{m|n}) \to U(\gl_{m|n})^{\otimes
l},  
\end{displaymath}
where $\Delta^{(l)}: Y(\gl_{m|n}) \to Y(\gl_{m|n})^{\otimes l}$ is the coproduct iterated $(l-1)$ times and
$\pi$ is the map given in \eqref{projection}. 
Then 
\begin{displaymath}
\kappa_l (t_{ij}^{(r)}) = 
\sum_{1 \le s_1 < \ldots < s_r \le l}\; \sum_{\substack{1 \le i_1, \ldots, i_{r-1} \le m+n}}
E_{i i_1}^{[s_1]} E_{i_1 i_2}^{[s_2]} \cdots E_{i_{r-1} j}^{[s_r]} 
(-1)^{\pa{i}\,+\, \pa{i}_1\, +\, \pa{i}_2 \,+\, \ldots + \pa{i}_{r-1}}
\end{displaymath}
where $E_{ij}^{[s]} = 1^{\otimes (s-1)} \otimes  E_{ij} \otimes 1^{\otimes (l-s)}$.  For any $r > l \ge 1$, we
have
$\kappa_{l}(t_{ij}^{(r)}) = 0$.
\begin{theorem}
Suppose we have fixed some ordering on the generators $t_{ij}^{(r)}$ $(1\le i,j \le m+n; \; r \ge 1)$ for the Yangian 
$Y(\gl_{m|n})$.   Then the ordered products of these, containing no second or higher order powers of the odd
generators, form a basis for $Y(\gl_{m|n})$.
\end{theorem}
\begin{proof}\!\!\footnote{This theorem was stated in \cite{Zhang} but the proof there is incomplete.} 
By relation
\eqref{dr4}, the graded algebra gr$_1 Y(\gl_{m|n})$ is
supercommutative, and thus the set of all ordered monomials in the generators $t_{ij}^{(r)}$ (with no second 
and higher order powers of the odd generators) span the Yangian $Y(\gl_{m|n})$.  It remains to show that they are
linearly independent.
We show that, for every $l \ge 1$,
the corresponding monomials in $\{ \kappa_l(t_{ij}^{(r)}) \; |\; 1\le r\le l\}$ are linearly independent in 
$\kappa_{l}(Y(\gl_{m|n}))$.  
Consider the filtration 
\begin{displaymath}
\mathrm{F}_0 U(\gl_{m|n})^{\otimes l} \subseteq \mathrm{F}_1 U(\gl_{m|n})^{\otimes l} 
\subseteq \mathrm{F}_2 U(\gl_{m|n})^{\otimes l} \subseteq \ldots
\end{displaymath}
on $U(\gl_{m|n})^{\otimes l}$ defined by setting each generator $E_{ij}^{[r]}$ to be of degree $1$. Then the
associated
graded algebra gr$\, U (\gl_{m|n})^{\otimes l}$ is the polynomial algebra on supersymmetric generators 
\[ x_{ij}^{[r]} := \mathrm{gr}_1 E_{ij}^{[r]}, \]
where  $x_{ij}^{[r]}$ is even if $\pa{i}+\pa{j} = \pa{0}$ and odd if $\pa{i} + \pa{j} = \pa{1}$.  
The map $\kappa_l$ preserves the filtration on the Yangian given by setting $\mathrm{deg}_1 (t_{ij}^{(r)})
= r$, and thus
defines a homomorphism between the corresponding graded algebras. 
It is enough to show that the same monomials in the elements
$y_{ij}^{(r)}~:=~\mathrm{gr}_r \kappa_{l} (t_{ij}^{(r)})$ 
in the graded algebra are linearly independent.  But for this, it is enough to show that the
superderivatives 
$d y_{ij}^{(r)}$ are linearly independent at a point.  We have:
\begin{equation*}
y_{ij}^{(r)} = \sum_{1 \le s_1 < \ldots < s_r < l}\; \sum_{1 \le i_1, \ldots, i_{r-1} \le n}
x_{i i_1}^{[s_1]} x_{i_1 i_2}^{[s_2]} \cdots x_{i_{r-1}j}^{[s_r]}
(-1)^{\pa{i} \,  + \pa{i}_1  + \ldots + \pa{i}_{r-1}}.
\end{equation*} 
We will show that the matrix $d \phi$ corresponding to the map 
$\left( d x_{ij}^{[s]}\right) \mapsto \left( d y_{ij}^{(r)} \right)$ has non-zero determinant at a point. 
It suffices to show that the determinant of this matrix is nonzero even when the variables are specialized
to $x_{kl}^{(s)} = \delta_{kl} c_s (-1)^{\pa{k}}$ for some distinct $c_s$ $(s \ge 1)$.  When the variables are
 specialized as described, we find:
\begin{equation*}
d y_{ij}^{(r)} = \sum_{s =1}^{l}\; \sum_{\substack{1 \le s_1 < \ldots < s_{r-1} \le l\\ s_i \ne s}}
c_{s_1} c_{s_2} \cdots c_{s_{r-1}} (-1)^{ \pa{i}} d x_{ij}^{[s]}.
\end{equation*}
Let $J$ be the
$(m+n) \times (m\times n)$ matrix 
$J=\left(\delta_{ij} (-1)^{\pa{i}}\right)$.  Then
$d\phi =  J\otimes X_l$,
where 
{\small \begin{displaymath}
X_l = \left( \begin{array}{cccc} 1& 1&  \ldots & 1\\
(c_2 + c_3 + \ldots +c_l)   & (c_1 + c_3 + \ldots + c_l)  & \ldots & (c_1 + c_2 + \ldots +c_{l-1}) \\
(\sum_{i,j \ne 1} c_{i} c_{j}) & (\sum_{i,j \ne 2} c_{i} c_{j})& \ldots &  (\sum_{i,j \ne l} c_{i}c_{j})\\
\vdots & \; & \; &\vdots \\
c_2c_3 \cdots c_l & c_1 c_3c_4 \cdots c_l & \ldots & c_2 c_3 \cdots c_{l-1}  
\end{array} \right). 
\end{displaymath}
}We show by induction that det$ X_l = \Pi_{1 \le i < j \le l} (c_i - c_j) \ne 0$, and hence
det$d\phi \ne 0$. Indeed, row-reducing $X_{l}$ gives the following matrix:
{\small \begin{displaymath}
\left( \begin{array}{cccc} 1 & \ldots & 1 & 1\\
(c_l -c_1)   & \ldots & (c_{l} - c_{l-1})      & 0\\
(c_l -c_1)\sum_{\substack{i,j \ne 1\\ i,j < l}} c_{i} c_{j} & 
\ldots &(c_{l} - c_{l-1})\sum_{\substack{i,j \ne l-1\\ i,j < l}} c_{i} c_{j}& 0 \\
\vdots & \; & \; & \vdots \\
(c_l - c_1)c_2 c_3 \cdots c_{l-1} & \ldots & (c_{l} - c_{l-1}) c_1 \cdots c_{l-2} & 0 
\end{array} \right),
\end{displaymath}
}which clearly has determinant $ (c_1 - c_l)(c_2 - c_l) \cdots (c_{l-1} - c_l ) $det$X_{l-1}$.  

Now, suppose we have some non-trivial linear combination $P$ of the ordered monomials in~$t_{ij}^{(r)}$ (with
no second or higher order powers of the odd generators) and take $l$ to be any number greater than
all the $r$ that occur in $P$.  Since the monomials in
$\kappa_{l}(t_{ij}^{(r)})$ are linearly independent in $\kappa_{l}(Y(\gl_{m|n}))$, we must have $\kappa_{l}(P)
\ne 0$. 
Therefore, $P \ne 0$ in the Yangian.
\end{proof} 

Now let $\gl_{m|n}[t]$ denote the algebra $\gl_{m|n} \otimes \mathbb{C}[t]$ with basis $\{E_{ij}t^r\}_{1\le i,j \le
m+n; r
\ge 0}$.  
\begin{corollary}\label{PBW}
The graded algebra gr$_2 Y(\gl_{m|n})$ is isomorphic to the algebra $U(\gl_{m|n}[x])$, via the map
\begin{eqnarray*}
\mathrm{gr}_2 Y(\gl_{m|n}) &\to&U(\gl_{m|n} [x]) \\
\mathrm{gr}_2^{r-1} t_{ij}^{(r)} &\mapsto& E_{ij} x^{r-1} (-1)^{\pa{i}} \quad (1 \le i,j \le m+n, r \ge 1).
\end{eqnarray*}
 \end{corollary}
\section{Gauss Decomposition of $T(u)$}\label{Gauss}
Here we describe a decomposition of the matrix $T(u)$ in
terms of the quasideterminants of Gelfand and Retakh \cite{GGRW}. 
\begin{defn}
Let $X$ be a square matrix over a ring with identity such that its inverse matrix $X^{-1}$ exists, 
and such that its $(j,i)$th entry is an invertible element of the ring.  Then the $(i,j)$th
\emph{quasideterminant} of $X$ is defined by the formula
{\small \begin{equation*}
|X|_{ij} = \left((X^{-1})_{ji}\right)^{-1} =: \left| \begin{array}{ccccc} x_{11} & \cdots & x_{1j} & \cdots & x_{1n}\\
&\cdots & & \cdots&\\
x_{i1} &\cdots &\boxed{x_{ij}} & \cdots & x_{in}\\
& \cdots& &\cdots & \\
x_{n1} & \cdots & x_{nj}& \cdots & x_{nn}
\end{array} \right|.
\end{equation*}
}\end{defn}
By Theorem 4.96 in
\cite{GGRW}, the matrix $T(u)$ defined in \eqref{T} has the following Gauss decomposition in terms of
quasideterminants:
\begin{equation*}
T(u) = F(u) D(u) E(u)
\end{equation*}
for unique matrices 
{\small
\begin{equation*}
D(u) = \left( \begin{array}{cccc} d_1 (u) & &\cdots & 0\\
& d_2 (u) &  &\vdots\ \\
\vdots & &\ddots &\\
0 &\cdots &  &d_{m+n} (u)
\end{array} \right),
\end{equation*}
\begin{equation*}
E(u)=\!\!\left( \begin{array}{cccc} \!\!1 &e_{12}(u) &\cdots & e_{1,m+n}(u)\\
&\ddots & &e_{2, m+n}(u) \!\!\!\!\\
& &\ddots & \vdots\\
0 & & &1 
\end{array} \right)\! ,\;
F(u) =\!\! \left( \begin{array}{cccc} \!\!1 & &\cdots &0\!\!\\
f_{21}(u) &\ddots & &\vdots\\
\vdots & & \ddots& \\
\!\!f_{m+n,1}(u) & f_{m+n, 2}(u) &\cdots &1\!
\end{array} \right)\!,
\end{equation*}
}where
{
\begin{eqnarray*}
d_i (u) &=& \left| \begin{array}{cccc} t_{11}(u) &\cdots &t_{1,i-1}(u) &t_{1i}(u) \\
\vdots &\ddots & &\vdots \\
t_{i1}(u) &\cdots &t_{i,i-1}(u) &\boxed{t_{ii}(u)}
\end{array} \right|, \\
e_{ij}(u) &=& d_i (u)^{-1} \left| \begin{array}{cccc} t_{11}(u) &\cdots &t_{1,i-1}(u) & t_{1j}(u) \\
\vdots &\ddots &\vdots & \vdots \\
t_{i-1,i}(u) &\cdots &t_{i-1,i-1}(u) & t_{i-1,j}(u)\\
t_{i1}(u) &\cdots &t_{i,i-1}(u) &\boxed{t_{ij}(u)}
\end{array} \right|,\\
f_{ji}(u) &=& \left| \begin{array}{cccc} t_{11}(u) &\cdots &t_{1, i-1}(u) & t_{1i}(u) \\
\vdots &\ddots &\vdots &\vdots \\
t_{i-1,1}(u) &\cdots &t_{i-1,i-1}(u) &t_{i-1,i}(u)\\
t_{ji}(u) &\cdots &t_{j, i-1}(u) &\boxed{t_{ji}(u)} 
\end{array} \right|
d_{i}(u)^{-1}.
\end{eqnarray*}
}
We use the
following notation for the coefficients:
\begin{eqnarray}
d_{i}(u) &=& \sum_{ r \ge 0} d_{i}^{(r)}u^{-r}; \quad
\left(d_{i}(u)\right)^{-1} = \sum_{r\ge 0} d_{i}^{_{\Large '}(r)}u^{-r};\\
e_{ij}(v)&=& \sum_{r \ge 1} e_{ij}^{(r)}v^{-r};\quad
f_{ji}(v)= \sum_{r \ge 1} f_{ji}^{(r)}v^{-r}.
\end{eqnarray}
It is easy to recover each generating series $t_{ij}(u)$ by multiplying together and taking
commutators of the series
 $d_i (u),\, e_j(u):= e_{j,j+1} (u)$, and  $f_j(v):=f_{j+1,j}(u)$
 for  $1 \le i \le m+n$,   $1 \le j \le m+n-1$.  Indeed,
 for each pair $i,j$ such that $1< i+1 <j \le m+n-1$, we have:
 \begin{equation}\label{egenerates}
 e_{ij}^{(r)} = (-1)^{\pa{j-1}} [e_{i,j-1}^{(r)}, e_{j-1}^{(1)}]; \quad 
 f_{ji}^{(r)} = (-1)^{\pa{j-1}} [f_{j-1}^{(1)}, f_{i,j-1}^{(r)} ] .
 \end{equation}
Thus the Yangian $Y(\mathfrak{gl}_{m|n})$ is generated by the coefficients of the series
\begin{equation*}
\left\{ d_{i}(u), e_{j}(u), f_{j}(u) \;|\; 1\le i \le m+n; 1\le j \le m+n-1\right\}.
\end{equation*} 
\section{Maps Between Yangians}
For Yangians $Y(\gl_{m|n})$ with small $m$ and $n$, such as $Y(\gl_{1|1})$ and $Y(\gl_{2|1})$, it is feasible
to use
this matrix relationship $T(u) = F(u) D(u) E(u)$ to translate the defining relations \eqref{dr2} into relations
between the generating series $d_{i}(u)$, $e_{j}(u)$ and $f_{j}(u)$.  However, in order to transfer these
results to the general case of $Y(\gl_{m|n})$ we must define various homomorphisms between Yangians.

\begin{lemma}The map $\rho_{m|n}: Y(\mathfrak{gl}_{m|n}) \to Y(\mathfrak{gl}_{n|m})$ defined by
\begin{equation*}
\rho_{m|n} (t_{ij}(u)) = t_{m+n+1-i, m+n+1-j} (-u).
\end{equation*}
is an associative algebra isomorphism.  
\end{lemma}
Note where we have swapped $m$ and $n$ in the above.  We use the same symbols for the generators of both
$Y(\gl_{m|n})$ and $Y(\gl_{n|m})$.  It should be clear from the context which algebra $t_{ij}(u)$ belongs
to.
\begin{proof} We check that the map $\rho_{m|n}$ preserves the defining relation \eqref{dr2}.  
\end{proof}
\begin{proposition}\label{zeta} Let $\zeta_{m|n}: Y(\gl_{m|n}) \to Y(\gl_{n|m})$ be the associative algebra isomorphism
 given by $\zeta_{m|n} = \rho_{m|n} \circ \omega_{m|n}$, where $\omega_{m|n}$ is the $Y(\gl_{m|n})$
automorphism given by
\[\omega_{m|n}: T(u) \mapsto T(-u)^{-1}.\]
That is,
\[\zeta_{m|n}: t_{ij}(u) \mapsto t'_{m+n+1-i, m+n+1-j}(u) .\]
Then:
\begin{equation}
\zeta_{m|n}: \left\{ \begin{array}{lll} d_{i}(u)& \mapsto & \left(d_{m+n-i+1}(u)\right)^{-1},\\
e_{k}(u)& \mapsto & - f_{m+n-k}(u) , \\
 f_{k}(u)& \mapsto & - e_{m+n-k}(u) , 
\end{array} \right.
\end{equation}
for $1\le i \le m+n$ and $ 1\le k \le m+n-1$.
\end{proposition}
\begin{proof}
We multiply out the matrix products $$T(u) = F(u) D(u) E(u)$$ and $$T(u)^{-1} = E(u)^{-1} D(u)^{-1}
F(u)^{-1}.$$
These show that for all $1\le i < j \le m+n$,
\begin{eqnarray*}
t_{ii}(u)&=& d_{i}(u) + \sum_{k < i} f_{ik}(u) d_{k}(u) e_{ki}(u),\\
t_{ij}(u)&=& d_{i}(u) e_{ij}(u) + \sum_{k< i} f_{ik}(u) d_{k}(u) e_{kj}(u),\\
t_{ji}(u)&=& f_{ji}(u) d_{i}(u) + \sum_{k <i} f_{jk}(u) d_{k}(u) e_{ki}(u),
\end{eqnarray*} and
\begin{eqnarray*}
t'_{ii}(u)&=& d_{i}(u)^{-1} + \sum_{k >i} e'_{ik}(u) d_{k}(u)^{-1} f'_{ki}(u),\\
t'_{ij}(u)&=& e'_{ij}(u) d_{j}(u)^{-1} + \sum_{k > j} e'_{ik}(u) d_{k}(u)^{-1} f'_{kj}(u),\\
t'_{ji}(u)&=& d_{j}(u)^{-1} f'_{ji}(u) + \sum_{k >j} e'_{jk}(u) d_{k}(u)^{-1} f'_{ki}(u),
\end{eqnarray*} where
\[ e'_{ij} (u) = \sum_{i=i_{0}< i_{1} < \ldots < i_{s} = j} 
(-1)^{s} e_{i_{0}i_{1}}(u) e_{i_{1}i_{2}}(u) \cdots e_{i_{s-1}i_{s}}(u)\]
and
\[ f'_{ji}(u) = \sum_{i=i_{0} < i_{1} < \ldots < i_{s} = j} 
(-1)^{s} f_{i_{s}i_{s-1}}(u) \cdots f_{i_{2}i_{1}}(u) f_{i_{1}i_{0}}(u).\]
Then immediately we have $\zeta_{m|n}(d_{1}(u)) = d_{m+n, m+n}(u)^{-1}$, 
$\zeta_{m|n}( e_{1j}(u)) =f'_{m+n, m+n+1-j}(u) $, and $\zeta_{m|n}(f_{j1}(u)) = e'_{m+n+1-j, m+n}(u)$.  By
induction on $i$, we derive:
\begin{eqnarray*}
\zeta_{m|n}(d_{i}(u)) &=& \left( d_{m+n+1-i}(u)\right)^{-1},\\
\zeta_{m|n}(e_{ij}(u)) &=& f'_{m+n+1-i, m+n+1-j} (u),\\
\zeta_{m|n}(f_{ji}(u)) &=& e'_{m+n+1-j, m+n+1-i}(u).
\end{eqnarray*}
The result stated in the proposition is the special case of this where $j=i+1$.
\end{proof}
When it is reasonable we will write simply $\zeta$ for the map $\zeta_{m|n}$.
The map $\zeta_{m|n}$ restricts to the isomorphism  $U(\gl_{m|n}) \to U(\gl_{n|m})$
defined by
\[
E_{ij} \mapsto E_{m+n+1-i, m+n+1-j}.
\]
It can be calculated explicitly (using induction and basic
properties of quasideterminants) for any $1\le i,j \le m+n$ to give the following result:
\begin{equation*}
\zeta(t_{m+n+1-i, m+n+1-j}^{(r)}) =
 \sum_{\substack{r_{1}+ \ldots +r_{p} = r\\ r_{1}, \ldots, r_{p} > 0}} (-1)^{p} \sum_{k_{1}, \ldots, k_{p-1}=1}^{m+n} 
 t_{ik_{1}}^{(r_{1})} t_{k_{1}k_{2}}^{(r_{2})}\ldots t_{k_{p-1}j}^{(r_{p-1})}. 
\end{equation*}
Also, $\zeta$ is not a Hopf algebra map between the two Yangians, but
instead has the following property.
\begin{proposition}
Let $\tau: Y(\gl_{n|m}) \otimes Y(\gl_{n|m}) \to Y(\gl_{n|m}) \otimes Y(\gl_{n|m}) $ be the map given by
\[ \tau(y_{1} \otimes y_{2}) = y_{2} \otimes y_{1}\, (-1)^{\pa{y}_{1}\pa{y}_{2}} \]
for all homogeneous elements $y_{1}, y_{2} \in Y (\gl_{n|m})$.  Then:
\[ (\zeta\otimes \zeta) \circ \Delta=  \tau \circ \Delta \circ \zeta.\]
\end{proposition}
\begin{proof}
Recall that 
\[\Delta : T(u) \mapsto T_{[1]}(u) T_{[2]}(u),\]
where following \cite{MNO} we write 
\begin{eqnarray*}
T_{[1]}(u) &= &\sum_{i,j=1}^{m+n} t_{ij}(u) \otimes 1 \otimes E_{ij} (-1)^{\pa{j} (\pa{i} +1)}, \\
 T_{[2]}(u) &=&  \sum_{i,j=1}^{m+n} 1 \otimes t_{ij}(u) \otimes E_{ij} (-1)^{\pa{j}(\pa{i}+1)}.
\end{eqnarray*}
Then since $\Delta$ is an algebra homomorphism and we must have that
\[
\Delta : T(u)^{-1} \mapsto T_{[2]}(u)^{-1}T_{[1]}(u)^{-1},
\]
which gives explicitly:
\[ \Delta (t_{ij}'(u)) = \sum_{k=1}^{m+n} t'_{kj}(u) \otimes t_{ik}'(u) (-1)^{(\pa{i}+\pa{k})(\pa{j}+
\pa{k})}.\]
It is easy to see that this coincides with
$ ((\zeta\otimes \zeta) \circ \tau \circ \Delta \circ \zeta ) \,(t_{ij}'(u))$.
\end{proof}
Finally, let $\varphi_{m|n}: Y(\mathfrak{gl}_{m|n}) \hookrightarrow Y(\mathfrak{gl}_{m+k|n})$ be the inclusion
which sends each
$t_{ij}^{(r)}\in Y(\mathfrak{gl}_{m|n})$ 
to the generator $t_{k+i, k+j}^{(r)}\in Y(\mathfrak{gl}_{m+k|n})$; and 
let $\psi_k: Y(\mathfrak{gl}_{m|n}) \to Y(\mathfrak{gl}_{m+k|n})$ be the injective homomorphism defined by
\begin{equation}\label{psi}
\psi_k = \omega_{m+k|n} \circ \varphi_{m|n} \circ \omega_{m|n}.
\end{equation}
Then, for any $1 \le i,j \le m+n$ (see Lemma 4.2 of \cite{BK}) we have:
\begin{equation*}
\psi_{k} (t_{ij}(u)) = \left| \begin{array}{cccc} t_{11}(u) &\cdots &t_{1k}(u) &t_{1, k+j}(u)\\
\vdots &\ddots &\vdots &\vdots \\
t_{k1}(u) &\cdots &t_{kk}(u) &t_{k, k+j}(u)\\
t_{k+i, 1}(u) &\cdots &t_{k+i,k}(u) &\boxed{t_{k+i, k+j}(u)}
\end{array} \right|.
\end{equation*}
As an immediate consequence we have the following lemma.
\begin{lemma}\label{corpsi}
For $k, l\ge 1$, we have
\begin{eqnarray*}
\psi_{k} (d_l(u))&=& d_{k+l} (u),\;\\
\psi_{k} (e_l (u))&=&e_{k+l}(u),\\
\psi_{k}(f_l (u))&=& f_{k+l}(u).
\end{eqnarray*}
\end{lemma}
Notice that the map $\psi_k$ sends $t'^{\, (r)}_{ij} \in Y(\mathfrak{gl}_{m|n})$ to the element $t'^{\, (r)}_{k+i, k+j}$ in $Y(\mathfrak{gl}_{m+k|n})$. 
Thus the subalgebra $\psi_k (Y(\mathfrak{gl}_{m|n}))$  is generated by the elements $\{ t'^{\, (r)}_{k+s, k+t}\}_{s,t =1}^n$.  Then, by
\eqref{useful},
 all elements of this subalgebra commute with those of the subalgebra generated by the elements $\{ t_{ij}^{(r)} \}_{i,j =1}^{k}$.
This implies in particular that for any $i,j \ge 1$, the quasideterminants $d_i (u)$ and $d_j (v)$ commute.
\section{Gauss Decomposition of $Y(\gl_{2|1})$}
We begin by defining a presentation of the Yangian $Y(\gl_{2|1})$ using the Gauss decomposition.  We
will then use this to give the more general result in the next section.  We use
the matrix relationship $T(u) = F(u) D(u) E(u)$ to convert the defining relations \eqref{dr2} into relations
between the generating series $d_{i}(u)$, $e_{j}(v)$ and $f_{j}(v)$.
Note that in the Yangian $Y(\gl_{1|1})$, and in the Yangian $Y(\gl_{2})$, we have the following:
{\small  \begin{eqnarray}
T(u)
 &=& \left( \begin{array}{ll} d_{1}(u) & d_{1}(u)\,e_{1}(u)\\
f_{1}(u)d_{1}(u) & f_1 (u)d_1 (u) e_1 (u) + d_2(u)
\end{array} \right)  \label{monodromy2}\\
T(v)^{-1}&=& \left( \begin{array}{ll} 
d_{1}(v)^{-1}\! + e_{1}(v) d_{2}(v)^{-1} f_{1}(v) & -e_{1}(v)\, d_{2} (v)^{-1}\!\!\!\\
- d_{2}(v)^{-1} f_{1}(v) & \phantom{-}d_2(v)^{-1}
\end{array} \right). \label{monodromyinverse2}
\end{eqnarray} }
whereas in the Yangian $Y(\gl_{2|1})$,
{\small
\begin{eqnarray*}
T(u)&=&\left( \begin{array}{lll}
d_1 (u) &d_1 (u) e_1 (u) &d_1 (u) e_{13} (u) \\
f_1 (u) d_1 (u) &f_1 (u) d_1 (u) e_1 (u) + d_2 (u) &f_1 (u) d_1 (u) e_{13} (u) + d_{2} (u) e_2 (u)\\
f_{31}(u)d_1 (u)&f_{31} (u) d_3 (u) e_1 (u) + f_2 (u) d_3 (u) &*
\end{array} \right),\\
T(v)^{-1}\!\!\!\!\!&=&\left( \begin{array}{ccc}
*  & * &  (e_1 (v) e_2 (v) - e_{13} (v)) d_3 (v)^{\!-1}\!\!\\
* &d_2 (v)^{\!-1}\!+ e_2 (v) d_3 (v)^{-1} f_2 (v)\! & - e_{2} (v) d_3 (v)^{-1}\\
\!d_{3}(v)^{\! -1} (f_2 (v) f_1 (v)\!-\! f_{31}(v))\! &-d_3 (v)^{-1} f_2 (v) & d_3 (v)^{-1}
\end{array} \right).
\end{eqnarray*}
}
These expressions for the entries of $T(u)$ allow us to derive the following relations.
\begin{lemma}\label{qrlemma}
We have the following identities in $Y(\gl_{2|1})$:
\begin{eqnarray}
(u-v){[}d_i(u), e_j(v){]} &=& 
\left\{ \begin{array}{rl} (\delta_{i,j} -\delta_{i,j+1})\,d_i(u) (e_j(v) - e_j(u)), & \text{ if }j=1;\\
                                (\delta_{i,j} +\delta_{i,j+1})\,d_i(u) (e_j(v) - e_j(u)) ,& \text{ if }j=2;
\end{array}\right. \label{qrlemma1a}\\
(u-v) [d_{i}(u), f_{j}(v)] &=& 
\left\{ \begin{array}{cl} -\, (\delta_{ij} - \delta_{i,j+1})(f_{j}(v) - f_{j}(u)) d_{i}(u), &\text{ if }j=1;\\
                                  -\, (\delta_{ij} + \delta_{i,j+1})(f_{j}(v) - f_{j}(u)) d_{i}(u), &\text{
if }j=2;
\end{array}\right. \nonumber                      \\
(u-v){[}e_{j}(u), f_{k}(v){]} &=& 
\left\{  \begin{array}{cl} \delta_{jk}\left(d_j(u)^{-1}d_{j+1}(u)- d_{j}(v)^{-1} d_{j+1}(v)\right), &\text{
if }j=1 ;\\
                                  -\,\delta_{jk} \left(d_j(u)^{-1}d_{j+1}(u)- d_{j}(v)^{-1} d_{j+1}(v)\right),
&\text{ if }j=2;
\end{array}\right. \nonumber \\
(u-v) [ e_j(u), e_j(v)] &=& \left\{ \begin{array}{cl} (e_j(v) - e_{j}(u))^2, &\text{ if  }j=1;\\
                                                                        0, & \text{ if }j=2; 
\end{array}\right. \label{qrlemma3a}\\
(u-v) [ f_{j}(u), f_{j}(v)]&=& \left\{ \begin{array}{cl} -\,(f_{j}(v) - f_{j}(u))^{2}, & \text{ if }j=1;\\
                                                                                      0, & \text{ if }j=2;
\end{array} \right. \nonumber \end{eqnarray} \begin{eqnarray}
(u-v) [e_{1} (u), e_2 (v)]&=& e_1(u)e_2(v) - e_1(v)e_2(v) -e_{13}(u) + e_{13}(v), \nonumber \\
(u-v)[f_{1}(u), f_{2}(v)]&=& -\, f_{2}(v)f_{1}(u) + f_{2}(v)f_{1}(v) +f_{31}(u)- f_{31}(v),\nonumber \\
{[{[e_{i}(u), e_{j}(v)]}, e_{j}(w)]}&+&{[{[e_{i}(u), e_{j}(w)]}, e_{j}(v)] }\; =
\; 0, \; \text{ if }|i-j| = 1;
\nonumber\\
{[{[f_{i}(u), f_{j}(v)]}, f_{j}(w)]}&+&{[{[f_{i}(u), f_{j}(w)]}, f_{j}(v)] }\;=\; 0, \;\text{ if } |i-j| =
1;\nonumber
\end{eqnarray}
where unless otherwise indicated the indices $i,j,k$ range over $i=1,2,3$ and $j,k=1,2$. 
\begin{proof}
We give a proof of just the first equation \eqref{qrlemma1a}, since the rest are proven similarly.  First,
note that by the remarks at the end of the previous section, $d_{3}(u) = \psi_{2}(d_{1}(u))$ commutes with 
$e_{1}(v) = t_{11}(v)^{-1}t_{12}(v)$.  Similarly, $e_{2}(v) = \psi_{1}(e_{1}(v))$ commutes with $d_{1}(u)$.
Now consider the quasideterminants $d_{1}(u), d_{2}(u)$ and $e_{1}(v)$ in the algebra 
$Y(\gl_{2})[[u^{-1},v^{-1}]]$.  Here, we have the matrices $T(u)$, $T(v)^{-1}$ as in \eqref{monodromy2} and 
\eqref{monodromyinverse2}.  By \eqref{useful}, 
\[
(u-v)[t_{11}(u), t'_{12}(v)] = t_{11}(u)t'_{12}(v) + t_{12}(u)t'_{22}(v),\]
but this is the same as 
\[(u-v)[d_{1}(u), -e_{1}(v)d_{2}(v)^{-1}] = -d_{1}(u)e_{1}(v)d_{2}(v)^{-1} +
d_{1}(u)e_{1}(u)d_{2}(v)^{-1}.\]
Cancelling $d_{2}(v)$ on the right gives the desired equation when $i=j=1$, but in $Y(\gl_{2})[[u^{-1},
v^{-1}]]$. 
We deduce the relation in $Y(\gl_{2|1}) [[u^{-1}, v^{-1}]]$ by following the natural inclusion 
$Y(\gl_{2})\!\hookrightarrow Y(\gl_{2|1})$ which sends generators in $Y(\gl_{2})$ to those of the same
name in $Y(\gl_{2|1})$.  

For the result when $i=2,j=1$, we consider the commutator $[t'_{22}(u), t_{12}(v)]$
in the algebra $Y(\gl_{2})[[u^{-1},v^{-1}]]$ and make the same deduction.
For the case $j=2$, we find the relations between $d_{1}(u)$, $d_{2}(u)$ and $e_{1}(v)$  in the
algebra $Y(\gl_{1|1})[[u^{-1},v^{-1}]]$, and map these into the algebra $Y(\gl_{2|1})[[u^{-1},v^{-1}]]$, by
following $\psi_{1}: Y(\gl_{1|1}) \to Y(\gl_{2|1})$.  
\end{proof}
\end{lemma}
\begin{theorem}\label{quasipresentation}
The algebra $Y(\mathfrak{gl}_{2|1})$ is generated by the even elements $d_{1}^{(r)}$, $d_{2}^{(r)}$, $d_{3}^{(r)}$, $d_{1}^{\,\prime\,(r)}$, $d_{2}^{\,\prime\,(r)}$,
 $d_{3}^{\,\prime\,(r)}$, $e_{1}^{(r)}$, $f_{1}^{(r)}$, and odd elements $e_{2}^{(r)}, f_{2}^{(r)}$, with $r
\ge 1$, subject only to the following relations:
\begin{eqnarray}
d_i^{(0)} &=& 1,\nonumber \\
\sum_{t=0}^r d_{i}^{(t)} \, d_{i}^{\,\prime\,(r-t)} &= &\delta_{r0},\nonumber \\
{[} d_{i}^{(r)}, d_{l}^{(s)} {]} &=& 0, \label{qr0} \\
{[}d_i^{(r)}, e_j^{(s)} {]} &=&
\left\{  \begin{array}{cl}(\delta_{ij} -\delta_{i,j+1})\sum_{t=1}^{r-1} d_i^{(t)}e_j^{(r+s-1-t)} , & \text{
if }j=1;\\
(\delta_{ij} +\delta_{i,j+1})\sum_{t=1}^{r-1} d_i^{(t)}e_j^{(r+s-1-t)}, &\text{ if }j=2; \end{array} \right.
\label{qr1}\\
{[} d_i^{(r)},f_j^{(s)}{]}&=&\left\{ \begin{array}{cl} 
-(\delta_{i,j}\!-\!\delta_{i,j+1})\sum_{t=1}^{r-1} f_j^{(r+s-1-t)}d_i^{(t)},& \text{ if }j=1;\\
-(\delta_{i,j}\!+\!\delta_{i,j+1})\sum_{t=1}^{r-1} f_j^{(r+s-1-t)}d_i^{(t)},& \text{ if } j=2;
\end{array}\right. \label{qr2}\\
{[}e_j^{(r)},f_k^{(s)}{]} &=&\left\{ \begin{array}{cl} 
-\,\delta_{jk} \sum_{t=0}^{r+s-1} d_{j}^{\,\prime\,(t)} d_{j+1}^{(r+s-1-t)}, & \text{ if }j=1;\\
\,\delta_{jk} \sum_{t=0}^{r+s-1} d_{j}^{\,\prime\,(t)} d_{j+1}^{(r+s-1-t)}, & \text{ if }j=2;
\end{array}\right. \label{qr3}\\
{[} e_1^{(r)}, e_1^{(s+1)} {]} -{[} e_1^{(r+1)}, e_1^{(s)} {]} &= &e_1^{(r)} e_1^{(s)} + e_1^{(s)} e_1
^{(r)} ,
\label{qr4}\\
{[} f_1^{(r+1)}, f_1^{(s)} {]} -{[} f_1^{(r)}, f_1^{(s+1)} {]} &=& f_1^{(r)} f_1^{(s)} + f_1^{(s)} f_1^{(r)}
, \label{qr5}
\end{eqnarray}
\begin{eqnarray}
{[} e_2^{(r)} ,e_2^{(s)} {]} \;=\; 0,\quad {[} f_2^{(r)} ,f_2^{(s)} {]} &=& 0,\nonumber \\
{[}e_1^{(r+1)}, e_2^{(s)}{]} - {[} e_1^{(r)}, e_2^{(s+1)}{]}& =& e_1^{(r)} e_2^{(s)},\nonumber\\
{[}f_1^{(r+1)}, f_2^{(s)}{]} - {[} f_1^{(r)}, f_2^{(s+1)}{]} &= &- f_2^{(s)} f_1^{(r)},\nonumber\\
{[}{[} e_1^{(r)}, e_2^{(s)}], e_2^{(t)}{]} \;+ \;{[}{[} e_1^{(r)}, e_2^{(t)}], e_2^{(s)} {]} &=&0,\label{qr9}\\
{[}{[} f_1^{(r)}, f_2^{(s)}], f_2^{(t)}{]} \;+ \;{[}{[} f_1^{(r)}, f_2^{(t)}], f_2^{(s)} {]} &=&0,\nonumber\\
{[}{[} e_2^{(r)}, e_1^{(s)}], e_1^{(t)}{]} \;+ \;{[}{[} e_2^{(r)}, e_1^{(t)}], e_1^{(s)} {]} &=&0,\label{qr11}\\
{[}{[} f_2^{(r)}, f_1^{(s)}], f_1^{(t)}{]} \;+ \;{[}{[} f_2^{(r)}, f_1^{(t)}], f_1^{(s)} {]} &=&0\nonumber 
\end{eqnarray}
for all $i,l=1,2,3$, $j,k=1,2$ and all $r,s,t \ge 1$.
\end{theorem}
\begin{remark}
Relations \eqref{qr4} and \eqref{qr5} are equivalent to the following relations:
\begin{eqnarray*}
{[}e_i^{(r)}, e_i^{(s)} {]} &=& \sum_{t=1}^{s-1} e_{i}^{(t)}e_i^{(r+s-1-t)} - \sum_{t=1}^{r-1} e_i^{(t)}e_i^{(r+s-1-t)}\\
{[}f_i^{(r)}, f_i^{(s)} {]}&=& \sum_{t=1}^{r-1} f_i^{(r+s-1-t)}f_i^{(t)} - \sum_{t=1}^{s-1} f_i^{(r+s-1-t)} f_i^{(t)}
\end{eqnarray*}
\end{remark}
\begin{proof}  
We follow the method given in the proof of Theorem 5.2 in \cite{BK}.
First, we show that the corresponding
coefficients of quasideterminants in the Yangian satisfy the relations given in the Theorem.
The first three relations are obvious from the fact that the $d_i(u)$'s commute and the
definition of the series $d_{i}'(u):= (d_{i}(u))^{-1}$.
The rest follow from the relations in Lemma~\eqref{qrlemma}.  We show the proof of only \eqref{qr1} and
\eqref{qr5}
since the rest are derived similarly.

Observe that for any formal series $g(u) = \sum_{r\ge 0} g^{(r)} u^{-r}$ we have the identity
\begin{displaymath}
\frac{g(u)-g(v)}{u-v} = - \; \sum_{r,s \ge 1} g^{(r+s-1)} u^{-r} v^{-s}.
\end{displaymath}
Then, by \eqref{qrlemma1a}, 
 \begin{eqnarray*}
 [d_i(u), e_j(v)] &=&
(\delta_{ij} - (-1)^{\delta_{i,2}} \delta_{i,j+1})
(\sum_{t\ge 1} d_i^{(t)}u^{-t})(\sum_{p,s \ge 1} e_j^{(p+s-1)} u^{-p}v^{-s}).
\end{eqnarray*}
Taking coefficients of $u^{-r}v^{-s}$ gives \eqref{qr1}. 

Now consider \eqref{qrlemma3a}. In the case where $j=1$, this expands out as follows:
{\small 
\begin{eqnarray*}
(u-v)[e_1 (u), e_1(v)] &=& \left( \sum_{r \ge 1} e_{1}^{(r)}u^{-r} - \sum_{s \ge 1} e_1^{(s)} v^{-s} \right)^2\\
&=& - \sum_{r,s\ge 1} e_{1}^{(r)}e_{1}^{(s)} u^{-r}v^{-s} 
- \sum_{r,s \ge 1} e_1^{(s)}e_1^{(r)} u^{-r}v^{-s}\\
&&\quad + \sum_{r,s \ge 1} e_{1}^{(r)}e_{1}^{(s)} u^{-r-s}
+ \sum_{r,s \ge 1} e_1^{(r)}e_1^{(s)} v^{-r-s}. 
\end{eqnarray*}} Taking coefficients of $u^{-r}v^{-s}$ on both sides gives the relation \eqref{qr5}.  \endgraf
Now let $\widehat{Y}$ be the algebra defined by the relations in Theorem 2. We have shown that there is an
associative algebra homomorphism $\widehat{Y} \to Y(\gl_{2|1})$ taking 
each generator in $\widehat{Y}$ to the quasideterminant coefficient of the same name in the Yangian.  
By \eqref{egenerates} these elements generate the Yangian, so this homomorphism is surjective.
We will now show that the algebra $\widehat{Y}$ is spanned as a  vector space by certain monomials, and that
the images of these monomials form a basis for the Yangian $Y(\gl_{2|1})$. 
It follows that the homomorphism is an isomorphism.\endgraf
Let $e_{13}^{(r)}$ and $f_{31}^{(r)}$ be the elements of $\widehat{Y}$ defined by
\begin{displaymath}
e_{13}^{(r)} = {[}e_{1}^{(r)}, e_2^{(1)}{]}, \quad f_{31}^{(r)} = {[f_{1}^{(r)}, f_{2}^{(1)}]} \quad
\text{ (c.f. \eqref{egenerates}). }
\end{displaymath}
We want to show that the algebra $\widehat{Y}$ is spanned by the set of ordered monomials in 
\begin{displaymath}
\{f_{31}^{(r)},f_{2}^{(r)},   f_{1}^{(r)},
 d_{1}^{(r)}, d_{2}^{(r)}, d_{3}^{(r)}, 
 e_{1}^{(r)}, e_{2}^{(r)}, e_{13}^{(r)}\, | \;r\ge1\},
\end{displaymath} taken in order some order so that that the $f$'s come before all the $d$'s, which come
before all the $e$'s.  It is clear from the relations \eqref{qr0}, \eqref{qr1}, \eqref{qr2} and \eqref{qr3}
that the
monomials in the above elements, where $f$'s come before $d$'s and $d$'s come before $e$'s, with the $d$'s
taken in some fixed order, do indeed span $\widehat{Y}$.

So our problem is to show that the subalgebra $\widehat{Y}^+$ of $\widehat{Y}$ generated by
elements $\{e_{i}^{(r)}\}_{i=1,2}$ is spanned by the monomials in $\{e_{1}^{(r)}, e_{2}^{(r)},
e_{13}^{(r)}; \; r\ge 1\}$ taken in some fixed order, and
similarly that the subalgebra $\widehat{Y}^{-}$ generated by elements $\{f_{i}^{(r)}\}_{i=1,2}$ is spanned
by the monomials in $\{f_{31}^{(r)},f_{2}^{(r)},   f_{1}^{(r)};\; r\ge 1\}$ taken in some fixed order. 
Consider $\widehat{Y}^+$. Define a filtration
\begin{displaymath}
L_0 \widehat{Y}^+ \subseteq L_1 \widehat{Y}^+ \subseteq \cdots
\end{displaymath}
on $\widehat{Y}^+$ by setting the degree of $e_{i}^{(r)}$ equal to $(r-1)$.  Let $gr^{L} \widehat{Y}^+$ 
be the associated graded algebra, and let $\pa{e_{i}}^{(r)}:= gr_{r-1}^{L} e_{i}^{(r)} \in gr^L \widehat{Y}^+$
for each $i= 1,2, 13$. Then we have the following:
\begin{eqnarray*}
{[}\pa{e}_{1}^{(r)}, \pa{e}_{1}^{(s)} {]} &=& 0,\quad
{[}\pa{e}_{2}^{(r)}, \pa{e}_{2}^{(s)} {]} \,=\, 0,\\
{[}\pa{e}_{13}^{(r)}, \pa{e}_{1}^{(s)} {]} &=& 0,\quad
{[}\pa{e}_{13}^{(r)}, \pa{e}_{2}^{(s)} {]} \,=\, 0,\\
{[}\pa{e}_{13}^{(r)}, \pa{e}_{13}^{(s)} {]} &=& 0,\quad
{[}\pa{e}_{1}^{(r)}, \pa{e}_{2}^{(s)} {]} \,=\, \pa{e}_{13}^{(r+s-1)}.
\end{eqnarray*}
Indeed, the first two identities are clear by the relations in the remark above. For the next two, first note
that
\begin{equation}\label{note}
[\pa{e}_1^{(r+1)}, \pa{e}_2^{(s)}] = [\pa{e}_1^{(r)}, e_2^{(s+1)}].
\end{equation}
Then 
\begin{eqnarray*}
{[}\pa{e}_{13}^{(r)}, \pa{e}_{12}^{(s)}{]} &=& {[}{[} \pa{e}_{12}^{(r)}, e_{23}^{(1)}{]}, e_{12}^{(s)}{]}
={[}{[} \pa{e}_{12}^{(1)}, e_{23}^{(r)}{]}, e_{12}^{(s)}{]}\\
&=& - {[}{[} e_{23}^{(r)} ,\pa{e}_{12}^{(1)}{]}, e_{12}^{(s)}{]}\;=\; - {[}{[} e_{23}^{(r)} ,\pa{e}_{12}^{(s)}{]}, e_{12}^{(1)}{]}\text{ (by \eqref{qr11})}\\
&=& - {[}{[} e_{23}^{(r+s-1)} ,\pa{e}_{12}^{(1)}{]}, e_{12}^{(1)}{]} \; =\;0 \text{ (by \eqref{qr11} again)}.
\end{eqnarray*}
Similarly,
\begin{eqnarray*}
{[} \pa{e}_{13}^{(r)}, \pa{e}_{23}^{(s)}{]} &=& {[}{[} \pa{e}_{12}^{(r)}, e_{23}^{(1)}{]}, e_{23}^{(s)}{]}\\
&=& -{[}{[} \pa{e}_{12}^{(r)}, e_{23}^{(s)}{]}, e_{23}^{(1)}{]} \text{ (by \eqref{qr9})}\\
&=& {[}{[} \pa{e}_{12}^{(r+s-1)}, e_{23}^{(1)}{]}, e_{23}^{(1)}{]} =0.
\end{eqnarray*}
The fifth relation is an easy consequence of these and the super-Jacobi identity:
\begin{equation*}
{[}\pa{e}_{13}^{(r)}, \pa{e}_{13}^{(s)}{]}
 \;=\; {[}{[} \pa{e}_{12}^{(r)}, \pa{e}_{23}^{(1)}{]}, \pa{e}_{13}^{(s)}{]}
= {[}{[} \pa{e}_{12}^{(r)}, \pa{e}_{13}^{(s)}, \pa{e}_{23}^{(1)}{]} {]} 
+ {[} \pa{e}_{12}^{(r)}, {[} \pa{e}_{23}^{(1)},\pa{e}_{13}^{(s)}{]}{]} =0
\end{equation*}
The final relation is just another extended application of \eqref{note}.
Given these calculations, it is clear that the graded algebra $gr^L \widehat{Y}^+$ is spanned by the set of all
ordered monomials in
$\{ \pa{e}_{ij}^{(r)} \}_{1\le i <j \le 3; r\ge 1}$ taken in some fixed order. Hence $\widehat{Y}^+$ is itself spanned by
the corresponding monomials in $\{e_{ij}^{(r)}\}_{1\le i <j \le 3; r\ge 1}$. The result for the subalgebra
$Y^{-}$ is shown similarly. \endgraf
Now we want to show that the monomials in
\begin{displaymath}
\{d_{i}^{(r)}\}_{1\le i\le3;\, r\ge 1} \cup \{ e_{ij}^{(r)}, f_{ji}^{(r)} \}_{1\le i<j \le 3; \, r\ge1}
\end{displaymath} taken in some fixed order so that $f$'s come before $d$'s and $d$'s come before $e$'s
form a basis for the Yangian $Y(\gl_{2|1})$.  By Corollary \ref{PBW}, we may identify the associated graded algebra 
gr$_2 Y(\gl_{m|n})$ with $U(\gl_{m|n}[t])$.  By the
definition of the quasideterminants, under this identification, gr$_2^{r} d_{i}^{(r+1)}$, 
gr$_2^{r} e_{ij}^{(r+1)}$, and gr$_2^{r} f_{ji}^{(r+1)}$ are identified, respectively, with $E_{ii} (-1)^{\pa{i}} t^r$, 
$E_{ij}(-1)^{\pa{i}} t^r$, and $E_{ji} (-1)^{\pa{j}} t^r$.  Then the result follows from the Poincar\'{e}-Birkhoff-Witt
theorem for Lie superalgebras (\cite{Scheunert}).
\end{proof}
\section{Gauss Decomposition of $Y(\gl_{m|n})$}
\begin{lemma}\label{mnlemma1}
The following relations hold in the algebra $Y(\gl_{m|n})[[u^{-1}, v^{-1}]]$.
\begin{eqnarray}
[d_{i}(u), d_{j}(v)]&=& 0 \text{ for all }1\le i,j\le m+n\\
(u-v){[d_{i}(u), e_{j}(v)]}&=& \left\{
\begin{array}{ll}
         (\delta_{ij} - \delta_{i,j+1}) d_{i}(u) (e_{j}(v) - e_{j}(u)),&{\text{ if }}\;j \le m-1,\\
         (\delta_{ij} + \delta_{i.j+1}) d_{i}(u) (e_{j}(v) - e_{j}(u)), &{\text{ if }}\;j = m,\\
        -(\delta_{ij} - \delta_{i,j+1}) d_{i}(u) (e_{j}(v) - e_{j}(u)),&{\text{ if }}\;j \ge m+1,
                                                                                         \end{array}\right.\\                                                                                
(u-v){[d_{i}(u), f_{j}(v)]}&=& \left\{
       \begin{array}{ll}
            - (\delta_{ij} - \delta_{i,j+1}) (f_{j}(v) - f_{j}(u)) d_{i}(u), &{\text{ if }}\; j \le m-1,\\
            - (\delta_{ij}+ \delta_{i,j+1}) (f_{j}(v) - f_{j}(u)) d_{i}(u), &{\text{ if }}\;j=m,  \\
              (\delta_{ij} - \delta_{i,j+1}) (f_{j}(v) - f_{j}(u)) d_{i}(u),& {\text{ if }}\;j \ge m+n-1,   
                                                                                            \end{array} \right. \\                                                                                                                                                
(u-v){[e_{i}(u), f_{j}(v)]}&=& 
         (-1)^{\pa{j+1}}\delta_{ij} \left(d_{i}(u)^{-1} d_{i+1}(u) - d_{i}(v)^{-1} d_{i+1}(v) \right),\\                                                                    
(u-v){[e_{j}(u), e_{j}(v)]}&=& \left\{
      \begin{array}{ll}  (-1)^{\pa{j+1}} \left( e_{j}(v) - e_{j}(u)\right)^{2},     &{\text{ if }}\;j \ne m,\\                                                                         
                                      0,                                        &{\text{ if }}\; j=m,   \end{array} \right. \\                                                                                                        
 (u-v){[f_{j}(u), f_{j}(v)]}&=& \left\{
     \begin{array}{ll}  -(-1)^{\pa{j+1}} \left( f_{j}(v) - f_{j}(u) \right)^{2},     &{\text{ if }}\; j \ne m,\\                                                                         
                                           0,                                     &{\text{ if }}\; j=m,          \end{array} \right.     \\
 (u-v){[e_{j}(u),e_{j+1}(v)]}\!&=& 
\! (-1)^{\pa{j+1}}\!\left( e_{j}(u)e_{j+1}(v)\! -\! e_{j}(v)e_{j+1}(v)\! -\!e_{j,j+2}(u) \!+\! e_{j,j+2}(v)
\right)\!,\;
\label{erelation}
\\                                                                               
(u-v){[f_{j}(u), f_{j+1}(v)]}\!&=& 
 \!\!-(-1)^{\pa{j+1}}\! \left( f_{j+1}(v)f_{j}(u) \!-\! f_{j+1}(v)f_{j}(v)\!-\!f_{j+2, j}(u)\!+\! f_{j+2,j}(v)\right)\!,\;\\
{[e_{i}(u), e_{j}(v)]}&=& 0 \text{ for }|i-j|>1;\\
{[f_{i}(u), f_{j}(v)]}&=&0 \text{ for }|i-j|>1;
\end{eqnarray}
\end{lemma}
\begin{proof}
The relations for $i,j$ between $1$ and $m$ are an easy consequence of those already found for the Yangians
$Y(\gl_{m})$ in \cite{BK} and $Y(\gl_{2|1})$ in Section \ref{qrlemma}, and the fact that the natural inclusions
$Y(\gl_{m})
\hookrightarrow Y(\gl_{m|n})$ and $Y(\gl_{2|1}) \hookrightarrow Y(\gl_{m|n})$ are homomorphisms. 
The remaining relations follow by applying the map $\zeta_{n|m}$ to the corresponding relations in $Y(\gl_{n|m})$.
\end{proof}
\begin{lemma}\label{mnlemma2} In addition, we have the following relations in $Y(\gl_{m|n})$ when $m >1$ and
$n>1$.  For any $r,s\ge 1$, 
\begin{equation}\label{extrarelation}
{[\,[e_{m-1}^{(r)}, e_{m}^{(1)}]\, , \, [e_{m}^{(1)}, e_{m+1}^{(s)}]\,]}\; = \; 0;\;\text{ and }\;\;
{[ \,[f_{m-1}^{(r)}, f_{m}^{(1)}]\, , \, [f_{m}^{(1)}, f_{m+1}^{(s)}]\, ]}\;=\; 0.
\end{equation}
\end{lemma}
\begin{proof}
We prove the result in $Y(\gl_{2|2})$, and then map this result
into the Yangian $Y(\gl_{m|n})$ via the map $\psi_{m-2}$.  
First we show the following relation:
\begin{equation}\label{preextra}
{[e_{13}(u)\,,\, e_{2}(z) e_{3}(z) - e_{2,4}(z)]}\; =\;0.                                                                                                                                                                                                                   
\end{equation}
Indeed, we have:
\begin{eqnarray*}
[e_{13}(u), \, e_{2}(z) e_{3}(z) - e_{24}(z)] &=&
 [e_{13}(u), e'_{24}(w)] \\
 &=& [t_{11}(u)^{-1} t_{13}(u)\,,\,- t'_{24}(w) t'_{44}(w)^{-1} ]\\
&=&0
\end{eqnarray*}
Now we find the commutator 
\[
(u-v)(w-z)[[e_{1}(u), e_{2}(v)]\, , \, [e_{2}(w), e_{3}(z)]].
\]
By \eqref{erelation}, this is
\begin{equation*}
{[e_{1}(u)e_{2}(v) - e_{1}(v)e_{2}(v) -e_{13}(u) + e_{13}(v) \;,\;
 -e_{2}(w)e_{3}(z) +e_{24}(w) +e_{2}(z)e_{3}(z) -e_{24}(z)]}\\.
 \end{equation*}
 Taking the coefficient of $u^{-r}z^{-s}$ and using \eqref{preextra} we find the first relation in
\eqref{extrarelation}.  The other part follows from this with the use of the map $\zeta$.
\end{proof}
Now we can state our main result.  The proof takes the same
line of reasoning as the proof of Theorem \ref{quasipresentation} but is somewhat longer and more
complicated. Again it is very closely based on the proof of Theorem 5.2 in \cite{BK}.
\begin{theorem}\label{presentation}
The Yangian $Y(\gl_{m|n})$ is isomorphic as an associative superalgebra to the algebra with
 even generators
$d_{i}^{(r)}$, $d_{i}^{\,\prime\,(r)}$, $f_{j}^{(r)}$, $e_{j}^{(r)}$, (for 
$1\le i \le m+n$, $1\le j \le m+n-1, j \ne m$, $r \ge1$) and odd generators $e_{m}^{(r)}$, $f_{m}^{(r)}$ (where again $r\ge1$) and the following defining relations:
\begin{eqnarray*}
d_i^{(0)} &=& 1; \\
\sum_{t=0}^r d_{i}^{(t)} \, d_{i}^{\,\prime\,(r-t)} &= &\delta_{r,0};\\
{[} d_{i}^{(r)}, d_{l}^{(s)} {]} &=& 0;
\end{eqnarray*}
\begin{eqnarray}
{[}d_{i}^{(r)}, e_{j}^{(s)}{]}& =& \left\{ \begin{array}{ll}
(\delta_{i,j}-\delta_{i,j+1})\sum_{t=0}^{r-1}d_{i}^{(t)}e_{j}^{(r+s-1-t)}, &\text{ for }\; 1\le j \le m-1,\\
(\delta_{i,j}+ \delta_{i,j+1})\sum_{t=0}^{r-1}d_{i}^{(t)}e_{j}^{(r+s-1-t)}, &\text{ for }\; j = m,\\
-(\delta_{i,j}- \delta_{i,j+1})\sum_{t=0}^{r-1}d_{i}^{(t)}e_{j}^{(r+s-1-t)}, &\text{ for }\; m+1\le j \le m+n-1,
\end{array}\right.
\label{r1}\\
{[d_{i}^{(r)}, f_{j}^{(r)}]} &=& \left\{ \begin{array}{ll}
-(\delta_{i,j}- \delta_{i,j+1})\sum_{t=0}^{r-1}f_{j}^{(r+s-1-t)}d_{i}^{(t)}, &\text{ for }\;1\le j\le m-1;\\
-(\delta_{i,j}+ \delta_{i,j+1})\sum_{t=0}^{r-1}f_{j}^{(r+s-1-t)}d_{i}^{(t)}, &\text{ for }\;j=m;\\
(\delta_{i,j}- \delta_{i,j+1})\sum_{t=0}^{r-1}f_{j}^{(r+s-1-t)}d_{i}^{(t)}, &\text{ for }\;m+1\le j\le m+n-1;
\end{array}\right.
\label{r2}\\
{[e_{j}^{(r)},f_{k}^{(s)}]}&=& \left\{ 
\begin{array}{ll}
  - \delta_{j,k} \sum_{t=0}^{r+s-1} d_{j}^{\,\prime\,(t)}d_{j+1}^{(r+s-1-t)}, & \text{ for }\; 1\le j\le m-1;\\
+\delta_{j,k} \sum_{t=0}^{r+s-1} d_{j}^{\,\prime\,(t)}d_{j+1}^{(r+s-1-t)}, & \text{ for }\; m\le j\le m+n-1;
\end{array}\right.   
\label{r3}
\end{eqnarray}
\begin{eqnarray}
{[e_{m}^{(r)}, e_{m}^{(s)}]} &= &0, \quad{[f_{m}^{(r)}, f_{m}^{(s)}]} \;= 0; \label{r4}      \\
{[e_{j}^{(r)}, e_{j}^{(s)}]} &=& (-1)^{\pa{j}} \left( \sum_{t=1}^{s-1} e_{j}^{(t)} e_{j}^{(r+s-1-t)} - 
                                 \sum_{t=1}^{r-1} e_{j}^{(t)} e_{j}^{(r+s-1-t)}\right), \text{ for }j\ne m;
 \label{r5}\\
{[f_{j}^{(r)}, f_{j}^{(s)}]} &=&  (-1)^{\pa{j}} \left( \sum_{t=1}^{r-1} f_{j}^{(t)} f_{j}^{(r+s-1-t)} - 
                               \sum_{t=1}^{s-1} f_{j}^{(t)} f_{j}^{(r+s-1-t)}\right), \text{ for }j \ne m;                              
\label{r5b}
\end{eqnarray}
\begin{eqnarray}
{[e_{j}^{(r)}, e_{j+1}^{(s+1)}]-[e_{j}^{(r+1)}, e_{j+1}^{(s)}]}&=& -(-1)^{\pa{j}} e_{j}^{(r)}e_{j+1}^{(s)}
\label{r7}\\
{[f_{j}^{(r+1)},f_{j+1}^{(s)}] - [f_{j}^{(r)}, f_{j+1}^{(s+1)}]}&=& -(-1)^{\pa{j}} f_{j+1}^{(s)}f_{j}^{(r)};
\label{r7b} \\
{[e_{j}^{(r)}, e_{k}^{(s)}]}\;=\;0;\; \text{ and }\; {[f_{j}^{(r)}, f_{k}^{(s)}]}&=&0, \text{ if } |j-k| >1;
\label{r9} \\     
{[}{[} e_j^{(r)}, e_k^{(s)}], e_k^{(t)}{]} \;+ \;{[}{[} e_j^{(r)}, e_k^{(t)}], e_k^{(s)} {]} \;&=&\;0, \text{ if }j \ne
k;\label{r10} \\
{[}{[} f_j^{(r)}, f_k^{(s)}], f_k^{(t)}{]} \;+ \;{[}{[} f_j^{(r)}, f_k^{(t)}], f_k^{(s)} {]} 
\;&=&\;0,\text{ if }j \ne k;\label{r10b}  \\
{[ \,[e_{m-1}^{(r)}, e_{m}^{(1)}]\, , \, [e_{m}^{(1)}, e_{m+1}^{(s)}]\, ]}&=& 0 \label{r11}\\
 {[ \,[f_{m-1}^{(r)}, f_{m}^{(1)}]\, , \, [f_{m}^{(1)}, f_{m+1}^{(s)}]\, ]}&=& 0 \label{r11b}
  \end{eqnarray}
for all $r,s,t \ge 1.$ and all admissible $i,j,k$. \end{theorem}
\begin{proof} Let $\widehat{Y}_{m|n}$ be the associative algebra given by the relations in the theorem.
By Lemma \ref{mnlemma1} and Lemma \ref{mnlemma2} the map from $\widehat{Y}_{m|n}$ to the Yangian 
$Y(\gl_{m|n})$ that sends every element of $\widehat{Y}_{m|n}$ to the element of the same name in the Yangian 
is a homomorphism.
We have already stated in Section \ref{Gauss} that $Y(\gl_{m|n})$ is generated by the elements:
\[\left\{ \left. d_{i}^{(r)}, e_{j}^{(r)}, f_{j}^{(r)} \right| 1\le i \le m+n, 1\le j \le m+n-1, r \ge 1 \right\}.
\]
Thus this homomorphism is surjective.  We need to show that it is injective. Our method is as follows: we show
that the algebra $\widehat{Y}_{m|n}$ is spanned as a vector space by the monomials in the elements 
$f_{ji}^{(r)}, d_{i}^{(r)}, e_{ij}^{(r)}$ with $1 \le i < j \le m+n$,\,$r\ge 1$, taken in some fixed order so
that the $f$'s come before $d$'s and $d$'s come before $e$'s.  (These elements are defined inductively by 
$f_{i+1,i}^{(r)} = f_{i}^{(r)};\quad e_{i,i+1}^{(r)} = e_{i}^{(r)}$ and 
\[f_{j,i}^{(r)} = [\, f_{j,j-1}^{(1)}\,,\, f_{j-1, i}^{(r)}\,]\, (-1)^{\pa{j-1}};\quad
 e_{i,j}^{(r)} = [\,e_{i,j-1}^{(r)}\,,\, e_{j-1,j}^{(1)}\,]\, (-1)^{\pa{j-1}},\quad 
\text{ for } j> i+1).
\]
Since the image of these monomials in the Yangian form a basis for $Y(\gl_{m|n})$, it follows that the map
is an isomorphism.

Let $\widehat{Y}^{+}_{m|n}$, $\widehat{Y}^{-}_{m|n}$ and $\widehat{Y}^{0}_{m|n}$
be the subalgebras of $\widehat{Y}_{m|n}$ 
generated by all elements of the form $e_{i}^{(r)}$, $f_{i}^{(r)}$ and $d_{i}^{(r)}$, respectively.  
By the defining relations \eqref{r1}, \eqref{r2} and  \eqref{r3}, we know that $\widehat{Y}_{m|n}$ is
spanned by the monomials where all $f$'s come before all $d$'s and all $d$'s come before all $e$'s.  Also,
since the $d$'s commute, we may assume that they are written in some fixed order.  If
we can show that the subalgebra $\widehat{Y}^{+}_{m|n}$ is spanned by the monomials in $e_{ij}^{(r)}$
written in some fixed order, then by applying the map $\zeta$ we can show that the subalgebra
$\widehat{Y}^{-}_{m|n}$ is similarly spanned by the monomials in $f_{ji}^{(r)}$ written in some fixed order.
This will then complete the proof.
 
Define an ascending filtration on $\widehat{Y}^{+}_{m|n}$ by setting deg$( e_{i}^{(r)}) = r-1$, and denote by
$gr^{L}\, \widehat{Y}^{+}_{m|n}$ the
corresponding graded algebra.  
Let $\pa{e}_{ij}^{(r)}$ be the image of $e_{ij}^{(r)}$ in the $(r-1)$-th
component of the graded algebra $gr^{L}\, \widehat{Y}^{+}_{m|n}$.  We claim that these images satisfy:
\begin{equation}\label{gradedcommutator}
{[\pa{e}_{ij}^{(r)}, \pa{e}_{kl}^{(s)}]} \; =
 \;(-1)^{\pa{j}} \, \delta_{kj} \, \pa{e}_{il}^{(r+s-1)} - (-1)^{\pa{i}\, \pa{j} + \pa{j}\, \pa{k} + \pa{i}\,\pa{k}}
\delta_{il}\, \pa{e}_{kj}^{(r+s-1)}. 
\end{equation}
From this relation it follows that the graded algebra $gr^{L}\, \widehat{Y}^{+}_{m|n}$ is spanned by the
monomials in $\pa{e}_{ij}^{(r)}$ taken in some fixed order.  Hence $\widehat{Y}^{+}_{m|n}$ is itself spanned
by the monomials in $e_{ij}^{(r)}$ taken in some fixed order.  

So now it remains only to prove the claim \eqref{gradedcommutator}. We begin by noting the following relations.
\begin{eqnarray}
{[\pa{e}_{i,i+1}^{(r)}, \pa{e}_{k,k+1}^{(s)}]}&=& 0, \text{ if } |i-k| \ne 1.\label{gr1}\\
{[\pa{e}_{i,i+1}^{(r+1)}, \pa{e}_{k,k+1}^{(s)} ]}&=& {[\pa{e}_{i,i+1}^{(r)}, \pa{e}_{k,k+1}^{(s+1)}]},
\text{ if } |i-k| = 1,\label{gr2}\\
{[\pa{e}_{i,i+1}^{(r)}\, ,\,  [\pa{e}_{i,i+1}^{(s)}, \, \pa{e}_{k,k+1}^{(t)}]\,]}&=& 
-{[\pa{e}_{i,i+1}^{(s)}\, ,\,  [\pa{e}_{i,i+1}^{(r)}, \, \pa{e}_{k,k+1}^{(t)}]\,]}, 
\text{ if }|i-k| = 1, \label{gr3}\\
\pa{e}_{ij}^{(r)}=\; {[\pa{e}_{i,j-1}^{(r)}, \, \pa{e}_{j-1,j}^{(1)}]} \, (-1)^{\pa{j-1}}
&=& {[\pa{e}_{i,i+1}^{(1)}, \, \pa{e}_{i+1,j}^{(r)}]}\, (-1)^{\pa{i+1}}, \text{ for }
j > i+1.
\label{gr4}
\end{eqnarray}
Here, \eqref{gr1} is a consequence of \eqref{r9}; \eqref{gr2} is a consequence of \eqref{r7}; and
\eqref{gr3} is a consequence of \eqref{r10}.  The first part of the last relation \eqref{gr4} follows from 
the definition of the elements $e_{ij}^{(r)}$.  The second part of \eqref{gr4} follows from
the first part using \eqref{gr2} and induction on the difference $j-i$.

Now we break up the problem of showing \eqref{gradedcommutator} into cases.  We assume without loss of
generality that $i\le k$.
If $j < k$, then $[\pa{e}_{ij}^{(r)}, \pa{e}_{kl}^{(s)}] = 0$ by \eqref{gr1} and \eqref{gr4}.
Consider the case where $j=k$.  By \eqref{gr2} and \eqref{gr4} we have 
\[ {[\pa{e}_{j-1, j}^{(r)}\, , e_{j, j+1}^{(s)}]} = (-1)^{\pa{j}} \pa{e}_{j-1, j+1}^{(r+s-1)}.\]
We bracket both sides of this with $\pa{e}_{j+1, j+2}^{(1)}$, $\pa{e}_{j+2, j+3}^{(1)}, \ldots$,
$\pa{e}_{l-1, l}^{(1)}$ in turn to obtain:
\[ {[\pa{e}_{j-1, j}^{(r)}, \, \pa{e}_{jl}^{(s)}]} \, =\, (-1)^{\pa{j}}\, \pa{e}_{j-1,l}^{(r+s-1)},\]
then bracket both sides of this new equation with $\pa{e}_{j-2, j-1}^{(1)}, \ldots,$ $\pa{e}_{i,i+1}^{(1)}$
to get the relation:
\[ {[\pa{e}_{i, j}^{(r)}, \pa{e}_{j, l}^{(s)}]} \,= \, (-1)^{\pa{j}}\, \pa{e}_{i, l}^{(r+s-1)}.\]

Before we consider the case $j>k$ in detail, we prove the following special cases:
\begin{eqnarray}
{[\pa{e}_{i, i+2}^{(r)}, \pa{e}_{i+1, i+2}^{(s)}]}\, =\, 0, \quad \text{ for } 1\le i \le m+n-2,\label{grlemma1}\\
{[\pa{e}_{i,i+1}^{(r)},\pa{e}_{i, i+2}^{(s)}]} \,=\,0, \quad \text{ for } 1\le i \le m+n-2,\label{grlemma1b}\\
{[\pa{e}_{i, i+2}^{(r)}, \pa{e}_{i+1, i+3}^{(s)}]} \,=\, 0\quad \text{ for } 1\le i\le m+n-3.\label{grlemma2}\\
{[\pa{e}_{ij}^{(r)}, \pa{e}_{k,k+1}^{(s)}]}\,=\, 0 \quad\text{ for }1\le i<k<j\le m+n \label{grlemma3}.
\end{eqnarray}  
Indeed, for \eqref{grlemma1}, we have:
\begin{eqnarray*}
\lefteqn{(-1)^{\pa{i+1}}\;{[\pa{e}_{i,i+2}^{(r)}, \pa{e}_{i+1,i+2}^{(s)}]} }\\
&=& {[[\pa{e}_{i,i+1}^{(r)}, \pa{e}_{i+1,i+2}^{(1)}], \pa{e}_{i+1,i+2}^{(s)}]} \; \text{ by \eqref{gr4}}\\
&=& -{ [[\pa{e}_{i,i+1}^{(r)},  \pa{e}_{i+1,i+2}^{(s)}],  \pa{e}_{i+1,i+2}^{(1)}]} \; \text{ by \eqref{gr3}}\\
&=&-{[[\pa{e}_{i,i+1}^{(r+s-1)}, \pa{e}_{i+1,i+2}^{(1)}], \pa{e}_{i+1,i+2}^{(1)} ]}  \; \text{ by
\eqref{gr2}},
\end{eqnarray*}
which is $0$ by \eqref{gr3}.  The relation \eqref{grlemma1b} is shown in a very
similar way.

When $i+1 = m,$ the relation \eqref{grlemma2} follows directly from \eqref{r11}. 
On the other hand, when $i+1\ne m$, the left-hand side of \eqref{grlemma2} equals
\begin{eqnarray*} 
\lefteqn{(-1)^{(\pa{i+1}\; +\; \pa{i+2})}\,{[[\pa{e}_{i,i+1}^{(r)}, \, \pa{e}_{i+1, i+2}^{(1)}]\,,\,
                               [\pa{e}_{i+1, i+2}^{(1)}, \pa{e}_{i+2, i+3}^{(s)}]\,]} }\\                              
 &=&  (-1)^{(\pa{i+1}\;+\; \pa{i+2})}
 \, [\pa{e}_{i+1, i+2}^{(1)},\, [\pa{e}_{i, i+1}^{(r)}, \pa{e}_{i+1, i+2}^{(1)}],\, \pa{e}_{i+2, i+3}^{(s)}]]\\
 &=&(-1)^{(\pa{i+1}\;+\; \pa{i+2})}
\,  {[\pa{e}_{i+1, i+2}^{(1)}, [\pa{e}_{i, i+1}^{(r)}, [\pa{e}_{i+1, i+2}^{(1)}, \pa{e}_{i+2,i+3}^{(s)}]]]}\\
&=& 
(-1)^{(\pa{i+1}\;+\; \pa{i+2})}
{[[\pa{e}_{i+1, i+2}^{(1)}, \pa{e}_{i, i+1}^{(r)}], [\pa{e}_{i+1, i+2}^{(1)}, \pa{e}_{i+2,i+3}^{(s)}]]},\\
&=& -\,\; (-1)^{(\pa{i+1}\;+\; \pa{i+2})}\,
{[[\pa{e}_{i, i+1}^{(r)},\pa{e}_{i+1, i+2}^{(1)}]\, ,\, [\pa{e}_{i+1, i+2}^{(1)}, \pa{e}_{i+2,i+3}^{(s)}]]}.
\end{eqnarray*}
Hence the commutator is zero.  Here we have used \eqref{r10} and the super-Jacobi identity, and the fact that
since $i+1 \ne m$, no two of the elements we are concerned with are odd.  

Finally, we use \eqref{gr4} relation to reduce the problem of showing \eqref{grlemma3} to that of showing
\begin{eqnarray*}
{[\pa{e}_{i, k+1}^{(r)}, \pa{e}_{k,k+2}^{(s)}]} &=& 0, \text{ and }\\
{[\pa{e}_{i, k+1}^{(r)}, \pa{e}_{k,k+1}^{(s)}]}&=&0,
\end{eqnarray*} 
for all $i\le k$.
The first of these relations follows from \eqref{grlemma1b} and \eqref{grlemma2} by induction on the difference
$k-i$, using
\eqref{gr4}.
The second follows from \eqref{grlemma1}, again by induction on $k-i$, using the relation \eqref{gr4}.

Now we properly begin the case $j>k$. We break this into the following subcases:
\begin{enumerate}
\item[Case 1:] $i<k$, $j=l$. 
Expanding $\pa{e}_{kj}^{(s)}$ by \eqref{gr4} and then using the super-Jacobi identity and \eqref{grlemma3}, we
have:
 \begin{equation*} 
 {[\pa{e}_{ij}^{(r)}, \, \pa{e}_{kj}^{(s)}]}=
 \pm [\,\pa{e}_{k,k+1}^{(1)}, [\pa{e}_{i,j}^{(r)},\pa{e}_{k+1, j}^{(s)}]\,].
  \end{equation*}
Continuing on in this fashion, we find:
\begin{equation*}
{[\pa{e}_{ij}^{(r)}, \, \pa{e}_{kj}^{(s)}]} \, =\, 
\pm \left[\pa{e}^{(1)}_{k,k+1}, \ldots[\pa{e}_{ij}^{(r)}, \pa{e}_{j-1, j}^{(s)}]\ldots\right],
\end{equation*}
so our problem reduces to showing that $[\pa{e}_{ij}^{(r)}, \pa{e}_{j-1, j}^{(s)}] = 0$.  We now expand out
the $\pa{e}_{ij}^{(r)}$ in this using \eqref{gr4} and apply the super-Jacobi identity to reduce this
problem
to that of showing that $[\pa{e}_{j-2,j}^{(r)}, \pa{e}_{j-1, j}^{(s)}] = 0$.  Then we have the result in this
case by \eqref{grlemma1}.
\item[Case 2:] $i<k$, $j>l$.  We expand out $\pa{e}_{kl}^{(s)}$ using \eqref{gr4} and then apply the
super-Jacobi identity and \eqref{grlemma3} to find:
\[ {[\pa{e}_{ij}^{(r)},\, \pa{e}_{kl}^{(s)}]} \,=\, \pm [ \pa{e}_{k,k+1}^{(1)}, \, [\pa{e}_{ij}^{(r)},
\pa{e}_{k+1, l}^{(s)}]].\]  Repeating this process as many times as is necessary we eventually get
\[ {[\pa{e}_{ij}^{(r)},\, \pa{e}_{kl}^{(s)}]} \,=\, \pm [ \pa{e}_{k,k+1}^{(1)}, \ldots,  [\pa{e}_{ij}^{(r)},
\pa{e}_{l-1, l}^{(s)}]\ldots ].\]  
which is $0$ by \eqref{grlemma3}.
\item[Case 3:] $i<k$, $j<l$.  We prove this case by induction on the difference $l-j$.  When $l-j=1$, we
have by expanding out $\pa{e}_{k,j+1}^{(s)}$ and using the super-Jacobi identity that
\begin{eqnarray*}
{[\pa{e}_{ij}^{(r)}, \, \pa{e}_{k,j+1}^{(s)}]}&=& 
[[\pa{e}_{ij}^{(r)}, \pa{e}_{kj}^{(s)}], \, \pa{e}_{j,j+1}^{(1)}] \,(-1)^{\pa{j}}
+ [\pa{e}_{kj}^{(r)}, \, [\pa{e}_{ij}^{(s)}, \, \pa{e}_{j,j+1}^{(1)}]]\, (-1)^{\pa{i}\, \pa{j} + \pa{j}\, \pa{k} +
\pa{i}\, \pa{k}}\\
&=& [[\pa{e}_{ij}^{(r)}, \pa{e}_{kj}^{(s)}], \, \pa{e}_{j,j+1}^{(1)}] \,(-1)^{\pa{j}}
+ [\pa{e}_{i,j+1}^{(s)}, \,\pa{e}_{kj}^{(r)}] (-1)^{(\pa{j}+\pa{j+1})(\pa{j}+ \pa{k})}.
\end{eqnarray*}
The first term is $0$ by the Case 1 and the second term is $0$ by Case 2.  When $l-j>1$, 
\begin{equation*}
{[\pa{e}_{ij}^{(r)}, \pa{e}_{kl}^{(s)}]} = 
[[\pa{e}_{ij}^{(r)}, \pa{e}_{k,l-1}^{(s)}], \, \pa{e}_{l-1, l}^{(1)}] (-1)^{\pa{l-1}},
\end{equation*} which is $0$ by the induction hypothesis.
\item[Case 4:] $i=k$, $j<l$.  We use \eqref{gr4} (and \eqref{gr1} and Case 2) to reduce this case to 
\eqref{grlemma1b}. 
\item[Case 5:] $i=k$, $j=l$. If $j=i+1$, then this is \eqref{gr1}.  Otherwise, we can expand out one term
with \eqref{gr4} to find:
\begin{equation*}
{[\pa{e}_{ij}^{(r)}, \, \pa{e}_{ij}^{(s)}]} = 
\pm [[\pa{e}_{i,j-1}^{(r)}, \pa{e}_{ij}^{(s)}],\, \pa{e}_{j-1,j}^{(1)}]
+ \pm [\pa{e}_{i,j-1}^{(r)}, \, [\pa{e}_{j-1, j}^{(1)}, \pa{e}_{ij}^{(s)}]\,].
\end{equation*}
The first term is $0$ by Case 4 and the second term is $0$ by Case 1.
\item[Case 6:] $i=k$, $j>l$.  This follows immediately from Case 4.
\end{enumerate}
This completes the proof of the claim \eqref{gradedcommutator}, which completes the proof of the theorem.
\end{proof}
\section{The Centre of $Y(\gl_{m|n})$}\label{centre}
The quantum Berezinian was defined by Nazarov \cite{Nazarov} as the following power series with coefficients in the Yangian $Y(\mathfrak{gl}_{m|n})$:
{\small 
\begin{eqnarray}
b_{m|n}(u)&:=&  \sum_{\rho \in S_m} \mathrm{sgn} (\rho) \,t_{\rho(1)1} (u) t_{\rho(2) 2} (u-1) \cdots
 t_{\rho(m)m} (u-m+1) \\
\nonumber && \times \;\sum_{\sigma \in S_n}\! \mathrm{sgn}(\sigma) \, t'_{m+1, m+\sigma(1)} (u-m+1) 
 \cdots t'_{m+n, m+\sigma(n)} (u-m+n) 
\end{eqnarray}
}
Recall from \cite{Gow} that we may also write the quantum Berezinian in the following form.
\begin{eqnarray}\label{qbear} 
b_{m|n}(u)
&=&  d_{1}(u)\, d_{2}(u-1) \cdots d_{m}(u-m+1) \\
\nonumber&& \times\, d_{m+1}(u-m+1)^{-1} \cdots d_{m+n} (u-m+n)^{-1}.
\end{eqnarray}
We shall prove that
the coefficients of this formal power series
generate the centre of the Yangian. This was conjectured by Nazarov who proved that the
quantum Berezinian was central \cite{Nazarov}.  A new proof of the centrality of the quantum Berezinian was
also given in \cite{Gow}.
\begin{lemma}\label{gradedcentre}
Let $\gl_{m|n}[x]$ be the polynomial current algebra and $I = E_{11} + \ldots + E_{m+n, m+n}$.  
The centre of $U(\gl_{m|n}[x])$ is generated by $I, Ix, Ix^{2}, \ldots$.
\end{lemma}
\begin{proof}
We reduce the problem to that of the well-known even case considered for example in Lemma 7.1 of \cite{BK}. 
First note that the supersymmetrization map gives an isomorphism between the
$\gl_{m|n}[x]$-modules $U(\gl_{m|n}[x])$ and $S(\gl_{m|n}[x])$, where $S(\gl_{m|n}[x])$ denotes the
supersymmetric algebra of $\gl_{m|n}[x]$.  The natural action of $\gl_{m|n}[x]$ on $S(\gl_{m|n}[x])$ is
obtained by extending the adjoint action.
The Lie algebra $\mathfrak{gl}_{m|n}$ has the root space decomposition:
\[
\mathfrak{gl}_{m|n} = \mathfrak{h} \oplus \bigoplus_{i=1}^{k} \mathfrak{g}_{\alpha_{i}}\,
\]
where $ \mathfrak{h}$ is the Cartan subalgebra, $\{\alpha_{1}, \ldots, \alpha_{k}\}$ is
the set of roots relative to $\mathfrak{h}$, and $\mathfrak{g}_{\alpha_{i}}$ is the one-dimensional root space
corresponding the root $\alpha_{i}$. 
Let
$e_{\alpha_{i}}$ be a root vector corresponding to root $\alpha_{i}$.
Suppose $P \in S(\gl_{m|n}[x])$ is an arbitrary $\gl_{m|n}$-invariant element and $M$ is the maximal integer
such that $e_{\alpha_{i}}x^{M}$ occurs in $P$ for some root $\alpha_{i}$.
Then we may write:
\begin{equation}
P= \sum_{\bf{s}} A_{\bf{s}} 
\left(  e_{\alpha_{1}} x^{M} \right)^{s_{1}} \ldots \left(  e_{\alpha_{k}} x^{M}\right)^{s_{k}},
\end{equation}   
where we sum over tuples of positive integers ${\bf{s}} = (s_{1}, \ldots, s_{k})$, and for each such
$\bf{s}$, the $A_{\bf{s}}$ is a monomial in 
elements $h x^{r}$ for $h \in \mathfrak{h}, r \ge 0,$ and $e_{\alpha_{i}}x^{r}$ for $r< M$.

For any $h\in \mathfrak{h}$, we have by assumption that:
\begin{eqnarray*}
0= {[hx, P]} &=& \sum_{{\bf{s}}} {[hx, A_{\bf s}]} 
\left(  e_{\alpha_{1}} x^{M} \right)^{s_{1}} \ldots \left(  e_{\alpha_{k}} x^{M}\right)^{s_{k}} \\
&& +\; \sum_{i=1}^{k} s_{i} \alpha_{i}(h)   \sum_{{\bf{s}}} A_{\bf{s}} 
\left(  e_{\alpha_{1}} x^{M} \right)^{s_{1}} \ldots \left(  e_{\alpha_{i}} x^{M} \right)^{s_{i}-1}
\ldots  \left(  e_{\alpha_{k}} x^{M}\right)^{s_{k}}\left( e_{\alpha_{i}}x^{M+1}\right).
\end{eqnarray*}
Then taking the coefficient of $\left( e_{\alpha_{i}}x^{M+1}\right)$ we find that for all $h\in
\mathfrak{h}$, and for all roots $\alpha_{i}$ that:
\[s_{i} \alpha_{i}(h)   \sum_{{\bf{s}}} A_{\bf{s}} 
\left(  e_{\alpha_{1}} x^{M} \right)^{s_{1}} \ldots \left(  e_{\alpha_{i}} x^{M} \right)^{s_{i}-1}
\ldots  \left(  e_{\alpha_{k}} x^{M}\right)^{s_{k}} = 0.\]
Since $\alpha_{i}(h)$ is not zero for all $h\in \mathfrak{h}$, and the monomials corresponding to different
${\bf s}$ are linearly independent, we must have that $s_{i} = 0$.  Thus $P$ is a sum of monomials in 
$hx^{r}$, where $h\in \mathfrak{h}$ and $r \ge 0$.  The Cartan subalgebra $\mathfrak{h}$ contains only even
elements, and so the action of $\mathfrak{gl}_{m|n}[x]$ on invariant elements $P$ is the same as the action of
$\gl_{m+n}[x]$. Then we may use Lemma 7.1 of \cite{BK} to obtain our desired result.
\end{proof}
\begin{theorem}\label{theoremcentre}
The coefficients of the quantum Berezinian generate the centre of $Y(\mathfrak{gl}_{m|n})$.
\end{theorem}
\begin{proof}
Write
\[ b_{m|n}(u) = 1+ \sum_{r\ge 1} b_{r} u^{-r}.\]
Our proof is based on that of Theorem 2.13. in \cite{MNO}.

Recall from Corollary \ref{PBW} that the graded algebra gr$_2 Y(\gl_{m|n})$ is isomorphic to $U(\gl_{m|n}[x])$. 
We show that for any $r =1,2,\ldots, $ the coefficient $b_{r}$ has degree $r-1$ with respect to deg$_{2} (.)$
and that its image in the $(r-1)$th component of gr$_{2} Y(\gl_{m|n})$ coincides with $I x^{r-1}$.  
Indeed, if we expand out the expression \eqref{qbear} for the quantum Berezinian, using the 
fact from \cite{GGRW} that 
\[
d_{j}(u) = 
t_{jj}(u) - \sum_{k,l < j} t_{jk}(u) (|T(u)_{\{1,2, \ldots , j-1\},\{1,2, \ldots , j-1\}}|_{lk})^{-1} t_{lj}(u),
\]
we find 
\[
b_{r} = \sum_{l_{1}+l_{2}+\ldots +l_{m+n} =r}
 t_{11}^{(l_{1})} t_{22}^{(l_{2})} \cdots t_{mm}^{(l_{m})}\cdot
 (-\,t_{m+1,m+1}^{(l_{m+1})}) \cdots (-\,t_{m+n, m+n}^{(l_{m+n})}) + \text{ terms of lower degree}.\]
Then it is clear that the terms  with $l_{i} = r$ for some $i=1,\ldots, m+n$ have degree $r-1$, and all else
have lower degree.  Then 
\[b_{r} = t_{11}^{(r)}+ \ldots + t_{mm}^{(r)} - t_{m+1, m+1}^{(r)} - \ldots - t_{m+n, m+n}^{(r)} + 
\text{ terms of lower degree}.
\]
The result follows when we evaluate the image of the graded part of this under the isomorphism in Corollary
\eqref{PBW}.
\end{proof}
\section{The Yangian $Y(\mathfrak{sl}_{m|n})$}\label{slmnsection}
Recall that the special linear Lie superalgebra $\,\mathfrak{sl}_{m|n}\,$ is the subalgebra of $\gl_{m|n}$
consisting of matrices with zero supertrace.
It may be defined explicitly by the 
following presentation \cite{GrozmanLeites, Scheunert2}.
We take generators $\{\, h_{i}, \; x^{+}_{j},\; x^{-}_{j} \, |\, 1\le i \le m+n-1 \}$.  The generators $\,h_{i},
\,
x_{j}^{\pm}\,$ are declared even for all $i$ and all $j \ne m$; the generators $\,{x_{m}^{\pm}}\,$ are declared
odd. The defining relations~are:
\begin{eqnarray*}
{[h_{i}, h_{j}] }&=& 0;\\
{[x_{i}^{+}, x_{j}^{-}]} &=& \delta_{i,j} h_{i};\\
{[h_{i}, x^{\pm}_{j}]} &=& \pm a_{ij} x^{\pm}_{j};\\
{[x_{m}^{\pm}, x_{m}^{\pm}]} &=& 0;\\
{[x_{i}^{\pm}, x_{j}^{\pm}]}&=& 0,\;\; \text{ if } |i-j| >1;\\
{[ x_{i}^{\pm}, \, [x_{i}^{\pm}, x_{j}^{\pm}]\,]} &=& 0,\;\;  \text{ if } |i-j| =1; \\
{[\,[x_{m-1}^{\pm}, x_{m}^{\pm}]\,,\, [x_{m+1}^{\pm}, x_{m}^{\pm}]]} &=& 0,
\end{eqnarray*}
for all $i,j$ between $1$ and $m+n-1$. Here $\,A \,=\, \left(a_{ij}\right)_{i,j= 1}^{m+n-1}\,$ is the \emph{symmetric} 
Cartan matrix of the Lie superalgebra $\mathfrak{sl}_{m|n}$, with entries $\;a_{ii} =2\,$ 
for all $i< m\,$; $\;a_{mm}=0\,$; $\;a_{ii} = -2\,$  for all $i > m\,$;
$\;\;a_{i+1,i} = a_{i,i+1} = -1\,$  for all $i<m\,$; $\;\;a_{i+1,i } = a_{i, i+1} = 1\,$  for all $i \ge m\,$;
 and all other entries are $0$. 

We define the Yangian $Y(\mathfrak{sl}_{m|n})$ associated to the special linear Lie superalgebra as the
following subalgebra of $Y(\gl_{m|n})$:
\[
Y(\mathfrak{sl}_{m|n}) : = \{ \; y \in Y(\gl_{m|n})\; |\; \mu_{f}(y) = y \text{ for all }f\;\},
\]
where we take $\mu_{f}$ as defined as in \cite{MNO}.  In other words, for a
formal power series 
\[
f = 1+ f_{1}u^{-1} + f_{2}u^{-2} +\ldots \quad \in \mathbb{C}[[u^{-1}]],
\]
the map $\mu_{f}$ is the automorphism of $Y(\gl_{m|n})$ given by
\[ \mu_{f}: T(u) \mapsto f(u) T(u).\]
This is justified by analogy with the definition of the Yangian $Y(\mathfrak{sl}_{N})$ as a subalgebra of
the Yangian $Y(\gl_{N})$ in \cite{MNO}.  Also, in the case where $m\ne n$ our definition
agrees with that arrived at by Stukopin \cite{StukopinA} through a quantization of the Lie bi-superalgebra 
$U(\mathfrak{sl}_{m|n}[x])$ (see Proposition~\ref{Stukopin}).  
\begin{proposition}
Let $Z_{m|n}$ denote the centre of the Yangian $Y(\gl_{m|n})$. Then for $m\ne n$, we have
\[ Y(\gl_{m|n}) \cong Z_{m|n} \otimes Y(\mathfrak{sl}_{m|n}).\]
\end{proposition}
\begin{proof}
We assume that $m>n$. (The result for $n<m$ follows from this by the application of the map $\zeta$).
The proof of this result is very similar to that of Proposition 2.16 in
\cite{MNO}. 
We use
the fact, stated there, that for any commutative associative algebra $\mathcal{A}$ and any formal series,
\[ a(u) = 1+ a_{1}u^{-1}+ a_{2}u^{-2}+ \ldots \in \mathcal{A}[[u^{-1}]],  \]
and any positive integer $K$ there exists a unique series
\[
\tilde{a}(u) = 1+ \tilde{a}_{1}u^{-1}+ \tilde{a}_{2}u^{-2}+ \ldots \in \mathcal{A}[[u^{-1}]]
\]
such that 
\begin{equation}\label{seriesexpansion}
 a(u) = \tilde{a}(u) \tilde{a}(u-1) \cdots \tilde{a}(u-K+1).
 \end{equation}
We take $a(u) = b_{m|n}(u)$ and $K= m-n$ in the commutative subalgebra  $Y^{0}\subset Y(\gl_{m|n})$
generated by the elements $d_{i}^{(r)}$ for $i=1,\ldots, m+n$ and $r\ge 1$.
Write 
\[ b_{m|n}(u) = \tilde{b}(u) \tilde{b}(u-1) \cdots \tilde{b}(u-m+n+1).\]
By the definition of the map $\mu_{f}$ we have that 
\[\mu_{f}(b_{m|n}(u)) = f(u)f(u-1)\cdots f(u-m+n+1) b_{m|n}(u).
\]
It follows from the uniqueness of the expansion \eqref{seriesexpansion} that $\mu_{f}(\tilde{b}(u)) = f(u)
\tilde{b}(u)$ for all $f$.  Also, the coefficients $\tilde{b}_{k},  (k\ge 1)$ of the series $\tilde{b}(u)$
generate the centre $Z_{m|n}$ since we may recover the coefficients of the series $b_{m|n}(u)$ from them. 
The remaining parts of the proof are exactly the same as in \cite{MNO}.
\end{proof}
\begin{lemma}\label{SYgenerators}
For any $m,n \ge 0$, the coefficients of the series
\begin{equation}
d_{1}(u)^{-1} d_{i+1}(u), \; e_{i}(u), \; f_{i}(u), \quad \text{ for }1\le i \le m+n-1,
\end{equation}
generate the subalgebra $Y(\mathfrak{sl}_{m|n})$.
\end{lemma}
\begin{proof}
It is clear that the coefficients of the series $d_{1}(u)$ together with those of the series listed above
generate the Yangian $Y(\gl_{m|n})$.  Also, for any $f$, the map $\mu_{f}$ leaves the coefficients of the
series in \eqref{SYgenerators} fixed and maps $\, \mu_{f}(d_{1}(u)) = f(u)d_{1}(u)$.  By the
Poincar\'e-Birkhoff-Witt theorem, any element $P$ of $Y(\gl_{m|n})$ is a polynomial in $\,d_{1}^{(1)},
d_{1}^{(2)}, d_{1}^{(3)}, \ldots\; $ and the other generators that are fixed by $\mu_{f}$ for all $f$.   
We can assume further that in each monomial in $P$ the generators are ordered so that the $f_{i}^{(r)}$'s come
before the $d_{i}^{(r)}$'s, which come before the $e_{i}^{(r)}$'s.  Suppose that $P \in
Y(\mathfrak{sl}_{m|n})$
and that $R$ is the maximum $r$ such
that $d_{1}^{(r)}$ occurs in $P$, and $K$ is the maximum power of $d_{1}^{(r)}$ occurring in $P$ for any $r$.
Fix $f= 1+\lambda u^{-R}$, where $\lambda$ is an arbitrary nonzero complex number.  
Then we can write:
\[ P = \sum_{a_{1}, a_{2}, \ldots, a_{R}} F_{a} D_{a} \left( d_{1}^{(1)}\right)^{a_{1}}
 \left(d_{1}^{(2)}\right)^{a_{2}} \ldots \left( d_{1}^{(R)}\right)^{a_{R}}
E_{a},\]
where $F_{a}, D_{a}$ and $E_{a}$ are monomials in the generators fixed by $\mu_{f}$, and we sum over all 
$R$-tuples $a= (a_{1}, a_{2}, \ldots, a_{R})$ of positive integers not exceeding $K$. 
Then
\[
\mu_{f}(P) =  \sum_{a_{1}, a_{2}, \ldots, a_{R}} F_{a} D_{a} \left(d_{1}^{(1)}\right)^{a_{1}}
\left(d_{1}^{(2)}\right)^{a_{2}} \ldots \left(d_{1}^{(R)}+\lambda \, \right)^{a_{R}}
E_{a} \;=\; P.
\]
By the linear independence of the different monomials and the fact that $\lambda$ is an arbitrary complex
number, we see that in fact $d_{1}^{(R)}$ cannot occur in $P
\in Y(\mathfrak{sl}_{m|n})$.
\end{proof}
Recall from  \cite{Kac77} that the family $A(m,n)$ of classical Lie superalgebras is defined by:
\begin{eqnarray}
 A(m-1,n-1) &=& \mathfrak{sl}_{m|n} \quad \text{ for } m \ne n;\quad  m,n \ge 1;\\
 A(n-1,n-1) &=& \mathfrak{sl}_{n|n}/ \langle I \rangle, \text{ for } n > 1,
\end{eqnarray}
where $\langle I \rangle$ is the one-dimensional ideal consisting of scalar matrices $\lambda I$, 
$(\lambda \in \mathbb{C})$. 
We define the Yangian of the classical Lie superalgebra $A(n-1,n-1)$ as the following quotient:
\begin{equation}
Y(A(n-1,n-1)) := Y(\mathfrak{sl}_{n|n})/ \left\langle b_{n|n} (u) =1 \right\rangle = Y(\mathfrak{sl}_{n|n})/B,
\end{equation}
where $B$ is the ideal in $Y(\mathfrak{sl}_{n|n})$ generated by the coefficients $b_{1}, b_{2}, \ldots$ of the
quantum Berezinian.  This definition is justified to a certain extent by Proposition \ref{justification}
below.
\begin{lemma}\label{centreofUsl}
For $n>1$, the centre of $U(A(n-1,n-1)[x])$ is trivial.  
\end{lemma}
\begin{proof}
We follow the argument of Lemma \ref{gradedcentre} using the properties of the
root-space decompostion given in \cite{Kac77}.
\end{proof}
\begin{proposition}\label{justification}
The centre of the Yangian $Y(A(n-1, n-1))$ is trivial.
\end{proposition}
\begin{proof}
We show that gr$_{2} Y(\mathfrak{sl}_{n|n}) \cong U(\mathfrak{sl}_{n|n}[x])$, and that 
\[ \text{gr} (Y(A(n-1,n-1))) \cong U(A(n-1,n-1)[x])).\]
Then the result follows from Lemma \ref{centreofUsl}.
Here we define the filtration on $Y(A(n-1, n-1))$,
\[ \mathbb{C} = A_{-1} \subset A_{0} \subset A_{1} \subset \ldots \subset A_{i} \subset \ldots, \]
by setting $A_{i} = Y_{i} + B$ where $Y_{i}$ is the set of
elements $ a \in Y(\mathfrak{sl}_{n|n})$ with deg$_{2}(a) \le i$, and gr$(Y(A(n-1,n-1)))$ is the
corresponding graded algebra.

The restriction of the map in Corollary \ref{PBW} to gr$_{2} Y(\mathfrak{sl}_{n|n})$ is injective onto its
image in 
$U(\mathfrak{sl}_{n|n})$.  By Lemma \ref{SYgenerators}, this is the image of the coefficients of the series 
$d_{1}(u)^{-1} d_{i+1}(u)$, $e_{i}(u)$ and $f_{i}(u)$, for $i =1,\ldots, 2n-1$. Now, for any $r\ge1$, the
coefficients of $u^{-r}$ these series are, respectively:
\begin{eqnarray*} 
&& t^{(r)}_{i+1, i+1} - t^{(r)}_{11} + \text{ elements of lower degree},\\
&& t^{(r)}_{i,i+1} + \text{ elements of lower degree},\\
&& t^{(r)}_{i+1, i} + \text{ elements of lower degree}.
\end{eqnarray*}
The image of these elements in $U(\gl_{n|n}[x])$ is:
\begin{eqnarray*}
&& (-1)^{\pa{i+1}} E_{i+1,i+1}x^{r-1} - E_{11}x^{r-1},\\ 
&& (-1)^{\pa{i}} E_{i,i+1} x^{r-1},\\
&& (-1)^{\pa{i+1}} E_{i+1, i} x^{r-1}.
\end{eqnarray*}
These elements generate precisely the subalgebra $U(\mathfrak{sl}_{n|n}[x])$.  
Thus we find that 
\[
{\text{gr}}_{2} Y(\mathfrak{sl}_{n|n}) \cong U(\mathfrak{sl}_{n|n}[x]).
\]
The natural projection map $p: Y(\mathfrak{sl}_{n|n}) \to Y(A(n-1,n-1))$ satisfies 
$p(Y_{i}) \subset A_{i}$,
and thus gives a natural surjective mapping 
\[ \text{gr}_{2} Y(\mathfrak{sl}_{n|n})\cong U(\mathfrak{sl}_{n|n}[x]) 
\to \text{gr} Y(A(n-1,n-1)),\]
with kernel the ideal ${\mathcal{I}} = \langle I, Ix, Ix^{2}, \ldots \rangle \subset U(\mathfrak{sl}_{n|n}[x])$.  Then
\[ \text{gr} Y(A(n-1, n-1)) \cong \, U(\mathfrak{sl}_{n|n}[x])/ {\mathcal{I}} \;  \cong U(A(n-1,n-1)[x]) .\]
\end{proof}
\begin{corollary}\label{centreofYslnn}
For $n>1$, the centre of the subalgebra $Y(\mathfrak{sl}_{n|n})$ is generated by the coefficients of the
quantum Berezinian $b_{n|n} (u)$.  
\end{corollary}
\section{Presentation of $Y(\mathfrak{sl}_{m|n})$}
Set \begin{eqnarray}\label{stukopingenerators1}
h_{i} ( u )  &=& d_i(u+ {\textstyle{\frac{1}{2}} }(-1)^{\pa{i}}(m-i)\;)^{-1} \nonumber
d_{i+1}(u+ {\textstyle{\frac{1}{2}} }(-1)^{\pa{i}}(m-i)\;),\\
x_{i}^{+} ( u ) & =& f_{i}(u+ {\textstyle{\frac{1}{2}} }(-1)^{\pa{i}}(m-i)\;)  \\
x_{i}^{-}(u ) &=&  (-1)^{\pa{i}}\,  e_{i}(u+ {\textstyle{\frac{1}{2}} }(-1)^{\pa{i}}(m-i)\;)\;  \nonumber
\end{eqnarray} for $1\le i\le m+n-1$, and use the following notation for the coefficients:
\begin{eqnarray}\label{stukopingenerators2}
h_i (u) &:= & 1+ \sum_{s\ge 0} h_{i,s} u^{-s-1},\quad \nonumber \\
x_{i}^{+} (u)&:=& \sum_{s\ge 0} x_{i,s}^{+} u^{-s-1},\\
x_{i}^{-} (u)&:= &\sum_{s\ge 0} x_{i,s}^{-} u^{-s-1}. \nonumber
\end{eqnarray}

Then we have the following presentation for the subalgebra $Y(\mathfrak{sl}_{m|n})$.
\begin{proposition}\label{Stukopin}
The subalgebra $Y(\mathfrak{sl}_{m|n})$ is isomorphic to the associative superalgebra over~{$\mathbb{C}$}
defined by the generators $x^{\pm}_{i,s}$ and $h_{i,s}$ for $1\le i \le m+n-1$ and
$s\in \mathbb{Z}_{+}$, and by the relations
\begin{eqnarray*}
{[}h_{i,r} ,h_{j,s}{]} &=& 0,\\
{[}x_{i,r}^{+} , x^{-}_{j,s}{]} &=& \delta_{ij} h_{i, r+s}\; ,\\
{[}h_{i,0} , x_{j,s}^{\pm}{]} &=&  \pm a_{ij} x_{j,s}^{\pm}\; ,\\
{[}h_{i,r+1}, x_{j,s}^{\pm} {]} - {[} h_{i,r},x_{j,s+1}^{\pm} {]} 
&=& \frac{\pm a_{ij}}{2}\; (h_{i,r}\, x_{j,s}^{\pm} + x_{j,s}^{\pm}\; h_{i,r}), \text{ for }i,j\text{ not
both }m, \nonumber\\
{[}h_{m,r+1}, x_{m, s}^{\pm} {]} &=& 0,\\
{[} x_{i, r+1}^{\pm}, x_{j,s}^{\pm} {]} - {[} x_{i,r}^{\pm} , x_{j, s+1}^{\pm} {]}
&=& \frac{\pm a_{ij}}{2} (x_{i,r}^{\pm} x_{j,s}^{\pm} + x_{j,s}^{\pm} x_{i,r}^{\pm} ),  \text{ for }i,j\text{ not
both }m,\nonumber \\
{[} x_{m, r}^{\pm}\,,\, x_{m,s}^{\pm} {]}&=&0,\\
{[ x_{i, r}^{\pm}\,,\, x_{j,s}^{\pm}{]}}&=& 0,\;\text{ if }|i-j| >1,\\
{[} x_{i,r}^{\pm} , {[} x_{i,s}^{\pm}, x_{j,t}^{\pm} {]}{]} + {[}x_{i,s}^{\pm} ,{[} x_{i,r}^{\pm} , x_{j,t}^{\pm}
{]}{]} &=& 0, \;\text{ if }|i-j|=1, \\
{[}\, {[} x_{m-1,r}^{\pm}, x_{m,0}^{\pm} {]}, {[}x_{m,0}^{\pm}  ,  x_{m+1, s}^{\pm}{]}\, {]} &=& 0,
\end{eqnarray*}
where $r$, $s$ and $t$ are arbitrary positive integers and $a_{ij}$ are the elements of the Cartan
matrix above.  The generators $x_{m, s}^{\pm}$ are odd and all other generators are even. 
\end{proposition}
\begin{proof}
For the duration of this proof we refer to the algebra given by the presentation in Proposition
\ref{Stukopin} as $\widetilde{Y}(\mathfrak{sl}_{m|n})$. By Lemma \ref{mnlemma1} we have a homomorphism 
$\varphi: \widetilde{Y}(\mathfrak{sl}_{m|n}) \to Y(\mathfrak{sl}_{m|n}) $ given by sending the elements 
$h_{i,s}, \, x_{i,s}^{\pm}$ to those defined in $Y(\mathfrak{sl}_{m|n})$ by \eqref{stukopingenerators1} and 
\eqref{stukopingenerators2}.  
By Lemma \ref{SYgenerators} this homomorphism is surjective.  We need to show $\varphi$ is injective. We
do this by constructing a set of monomials that span $\tilde{Y}(\mathfrak{sl}_{m|n})$, and whose image under
$\varphi$ is a basis for the Yangian $Y(\mathfrak{sl}_{m|n})$.  Following \cite{LevendorskiiPBW, StukopinA} we
construct this basis as follows.

Let $\alpha$ be a positive root of $\mathfrak{sl}_{m|n}$ and $\alpha = \alpha_{i_1}+ \ldots
+\alpha_{i_{p}}$ a decomposition of $\alpha$ into a sum of roots such that 
\[x_{\alpha}^{\pm} = [x_{i_{1}}^{\pm}, [x_{i_{2}}^{\pm}, \ldots, [x_{i_{p-1}}^{\pm},x_{i_{p}}^{\pm}]\ldots]]\]
is a nonzero root vector in $\mathfrak{sl}_{m|n}$. Suppose $s>0$ and we have a decomposition
$s= s_{1}+\ldots + s_{p}$ of $s$ into $p$ non-negative integers.  Then define the \emph{root vector}
$x_{\alpha, s_{1}+\ldots+ s_{p}}^{\pm}$ in the Yangian by
\begin{equation}\label{rootvectors}
x_{\alpha, s_{1}+\ldots + s_{p}}^{\pm} 
= [x_{i_{1}, s_{1}}^{\pm}, [x_{i_{2}, s_{2}}^{\pm}, \ldots,
 [x^{\pm}_{i_{p-1}, s_{p-1}},x_{i_{p}, s_{p}}^{\pm}]\ldots]].
\end{equation}

With respect to the second filtration defined in \eqref{filtrations}, the
degree of an element $h_{i,s}$ or $x_{i,s}^{\pm}$ is equal to its second index $s$, and 
deg$_{2} (x_{\alpha, s_{1}+\ldots + s_{p}}^{\pm}) = s$.  
If $s= s_{1}'+\ldots +s_{p}'$ is another decomposition of $s$ into non-negative integers, then (since
the defining relations in Proposition~\ref{Stukopin} are satisfied by the elements of the Yangian) we have 
\begin{equation}
\text{deg}_{2}(x_{\alpha, s_{1}'+\ldots +s_{p}'}^{\pm} - x_{\alpha, s_{1}+\ldots + s_{p}}^{\pm})\le s-1.
\quad 
\end{equation}
Now for each $s>0$ fix the decomposition $s= 0+\ldots +0+s$ to be used always and write 
$x_{\alpha, s}^{\pm} =x_{\alpha, 0+\ldots + 0+s}^{\pm}$.
Also any positive root  $\alpha$ is just $\alpha = \epsilon_{i} - \epsilon_{j}$ for some $1\le i \le
j-1 \le m+n-1$.  We then write: $x_{i,j; s}^{\pm} =x_{\alpha, 0+\ldots + 0+s}^{\pm}$.
Now choose
any total ordering $\prec$ on the set 
\[
\{x_{i,j; q}^{-}, h_{i, r}, x_{i,j; s}^{+}\; |\;1\le i\le j-1 \le m+n-1,  \; q,r,s>0
\}\]
and define $\Omega (\prec)$ to be the set of ordered monomials in these elements, where the odd elements 
($x_{i,j;r}^{\pm}$ with $i\le m$ but $j>m$) occur with power at most $1$.  

Define the length $l(M)$ of a monomial in $x_{i,j; q}^{-}, h_{i, r}, x_{i,j; s}^{+}$ as the number of
factors of $M$ and note that by the relations in Proposition~\ref{Stukopin}, if we rearrange the factors of
$M$, then
we obtain additional terms of either smaller degree, or the same degree but smaller length.  Then by induction
on the degree $d$ of a polynomial, and for fixed degree $d$, induction on the maximal length of its terms,
we see that $Y(\mathfrak{sl}_{m|n})$ is spanned by the elements of $\Omega (\prec)$.  (This argument is
given in \cite{LevendorskiiPBW} for the Yangian $Y(\mathfrak{sl}_{N})$).

Now suppose that some linear combination $\Sigma$ of the monomials in $\Omega (\prec)$ is equal to $0$, and
that the highest degree of a monomial term in $\Sigma$ is $r$. 
The degree $r$ part of $\Sigma$ must be equal to zero.  This will be the sum of 
products of the highest degree parts of elements $x_{i,j; r}^{-}, h_{i, r}, x_{i,j; r}^{+}$, which by 
the isomorphism gr$_{2}Y(\mathfrak{sl}_{m|n}) \cong U(\mathfrak{sl}_{m|n}[x])$ get mapped to the elements 
\[
\varepsilon_{i,j}^{-}\,E_{ij} x^{r-1}, \quad
\; (-1)^{\pa{i}}E_{ii} - (-1)^{\pa{i+1}} E_{i+1, i+1}; \quad
\varepsilon_{i,j}^{+} \, E_{ji} x^{r},
\]
respectively, where $\varepsilon_{i,j}^{\pm}$ is some power of $-1$.  
Together these elements form basis for
$\mathfrak{sl}_{m|n}[x]$, and so by the PBW theorem
for Lie superalgebras (\cite{Scheunert}) the set of ordered monomials in these, containing powers of at most
one of the odd elements, are linearly independent.  This implies that the highest degree part of $\Sigma$
must in fact be trivial.  Thus $\Omega (\prec)$ is a basis for $Y(\mathfrak{sl}_{m|n})$.

Now, we define a set $\widetilde{\Omega}(\prec)$ in $\widetilde{Y}(\mathfrak{sl}_{m|n})$ by the same
formulas as in \eqref{rootvectors}, except now we take the symbols to represent the elements of 
$\widetilde{Y}(\mathfrak{sl}_{m|n})$.  We define a filtration on $\widetilde{Y}(\mathfrak{sl}_{m|n})$ by
setting the degree of an element $h_{i,s}$ or $x_{i,s}^{\pm}$ equal to its second index $s$.  All the
arguments required to show that $\Omega(\prec)$ span the Yangian depended only on the relations in
Proposition~\ref{Stukopin}, and thus hold true for $\widetilde{\Omega}(\prec)$ in
$\widetilde{Y}(\mathfrak{sl}_{m|n})$.  Then $\widetilde{\Omega}(\prec)$ is a set of monomials that span
$\widetilde{Y}(\mathfrak{sl}_{m|n})$,  and whose image under $\phi$, $\Omega(\prec)$, is a basis for
$Y(\mathfrak{sl}_{m|n})$. 
\end{proof}
This is the presentation given by Stukopin \cite{StukopinA, Stukopin}, except that the last
relation has been corrected. Stukopin derives this presentation of the Yangian $Y(\mathfrak{sl}_{m|n})$
according to the definition of Yangian given in \cite{Drinfeld3}, as the quantization of the Lie 
bi-superalgebra $\mathfrak{sl}_{m|n}[t]$.   He names it after the series of classical Lie
superalgebras $A(m-1,n-1)$ and defines it only for the case $m \ne n$, since in the case where $m=n$ the Lie
superalgebra $\mathfrak{sl}_{m|n}[x]$ does not have a canonical Lie bi-superalgebra structure.  Stukopin
defines the root vectors given in the proof of Proposition~\ref{Stukopin} and 
gives a Poincar\'e-Birkhoff-Witt theorem for the Yangian $Y(\mathfrak{sl}_{m|n})$ using the same general 
argument as Levendorskii \cite{
LevendorskiiPBW}. The linear independence part of this PBW theorem may now also be obtained as a corollary of 
Proposition~\ref{Stukopin}.
%


\subsection*{Acknowledgements}
Thanks to my PhD supervisor Alex Molev for his patient assistance. 
He made many 
suggestions for improving this paper.  Thanks also to Vladimir Stukopin for explaining details of his work
to me through email.
\bibliographystyle{plain}
\def\cprime{$'$} \def\cprime{$'$} \def\cprime{$'$} \def\cprime{$'$}
  \def\cprime{$'$} \def\cprime{$'$} \def\cprime{$'$} \def\cprime{$'$}
\end {document}